\documentclass[12pt]{amsart}
\usepackage{amssymb,amscd}

\newtheorem{theorem}{Theorem}
\newtheorem{proposition}[theorem]{Proposition}

\newtheorem{problem}[theorem]{Problem}
\newtheorem{lemma}[theorem]{Lemma}
\newtheorem{corollary}[theorem]{Corollary}

\theoremstyle{remark}
\newtheorem{remark}[theorem]{Remark}
\newtheorem{example}[theorem]{Example}

\theoremstyle{definition}
\newtheorem{definition}[theorem]{Definition}
\newtheorem*{convention}{Convention}

\numberwithin{equation}{section}
\numberwithin{theorem}{section}

\newcommand{\A}{\mathcal{A}}
\newcommand{\B}{\mathcal{B}}
\newcommand{\C}{\mathbb{C}}

\newcommand{\e}{\varepsilon}
\newcommand{\G}{\mathcal{G}}
\newcommand{\h}{\mathfrak{H}}

\newcommand{\kom}{\mathcal{K}}
\newcommand{\kH}{\mathfrak{K}}

\newcommand{\M}{\mathcal{M}}
\newcommand{\NN}{\mathcal{N}}

\newcommand{\ozm}{\Omega(\z(\M))}
\newcommand{\Os}{\mathcal{O}}
\newcommand{\p}{\mathcal{P}}
\newcommand{\ps}{(\p(\M)/\sim~)}
\newcommand{\Q}{\mathbb{Q}}
\newcommand{\RR}{\mathcal{R}}
\newcommand{\R}{\mathbb{R}}
\newcommand{\U}{\mathcal{U}}
\newcommand{\z}{\mathcal{Z}}

\begin{document}

\title[Unitary orbits]{Unitary orbits of normal operators in von Neumann algebras}
\author{David Sherman}
\address{Department of Mathematics\\ University of California\\ Santa Barbara, CA 93106}
\email{dsherman@math.ucsb.edu}

\begin{abstract}
The starting points for this paper are simple descriptions of the norm and strong* closures of the unitary orbit of a normal operator in a von Neumann algebra.  The statements are in terms of spectral data and do not depend on the type or cardinality of the algebra.

We relate this to several known results and derive some consequences, of which we list a few here.  Exactly when the ambient von Neumann algebra is a direct sum of $\sigma$-finite algebras, any two normal operators have the same norm-closed unitary orbit if and only if they have the same strong*-closed unitary orbit if and only if they have the same strong-closed unitary orbit.  But these three closures generally differ from each other and from the unclosed unitary orbit, and we characterize when equality holds between any two of these four sets.  We also show that in a properly infinite von Neumann algebra, the strong-closed unitary orbit of any operator, not necessarily normal, meets the center in the (non-void) left essential central spectrum of Halpern.  One corollary is a ``type III Weyl-von Neumann-Berg theorem" involving containment of essential central spectra.
\end{abstract}

\maketitle

\section{Introduction} \label{S:intro}

Any research involving unitary orbits of Hilbert space operators is necessarily related to an enormous amount of mathematical literature.  The passage from operator to unitary orbit is a natural projectivization: by putting all equivalent representations on an equal footing, one isolates the algebraic and measure-theoretic information present in the operator.  If one considers the norm-closed unitary orbit instead, the useful relation of approximate equivalence comes into play.  Larger closures are related to homomorphic images and operator models.  There are deep generalizations to equivalence of representations and approximate innerness of completely positive maps.

The analogous passage from operator to unitary orbit inside a general von Neumann algebra $\M$ is no less natural.  When $\M$ is represented on a Hilbert space, one may view the difference between this procedure and the preceding one as the preservation of certain extra symmetries which can be realized as the operators in $\M'$.  But our viewpoint in this paper is almost entirely \textit{non-spatial}: the universe is $\M$, and no Hilbert space is needed.  The main goal is to give and apply simple descriptions of the norm and strong* closures of the unitary orbit of a normal operator in a von Neumann algebra.

Let us introduce a bit of notation: in any von Neumann algebra $\M$, we write $\U(\M)$ for the unitary group and $\U(h)$ for ${\{uhu^* \mid u \in \U(\M) \}}$, the unitary orbit of $h \in \M$.  We also write $\chi_E(h)$ for the spectral projection of normal $h$ corresponding to $E \subseteq \C$. 

Of course the first natural question to ask is ``When do two normal operators $h$ and $k$ have the same unitary orbit?"  For $\B(\h)$, the von Neumann algebra of all bounded operators on a Hilbert space $\h$, this was a motivating force for the development of spectral theory through the first half of the twentieth century.  Names often mentioned in this context include Hellinger, Hahn, Wecken, Plessner, Rokhlin, Nakano, and Halmos.  The machinery they developed, \textit{multiplicity theory}, completely solves this problem and generalizes to much broader situations, but it is fundamentally inadequate for classifying unitary orbits in von Neumann algebras.  Further comments about this question - mostly an outline of the difficulties - are postponed until Section \ref{S:agree}.
 
It turns out that the classification of closed unitary orbits is much more tractable.  This problem also has its origins in the first half of the twentieth century, in work on the norm and strong continuity of the map that sends a self-adjoint Hilbert space operator to its spectral resolution.  See the theorems and references in \cite[Section 135]{RSN}.  But precise descriptions of closed unitary orbits of normal Hilbert space operators really date from the 1970s.      

\begin{theorem} \label{T:normbh} $($\cite{GP,H1974,H1981}$)$
The following conditions are equivalent for normal operators $h$ and $k$ in $\B(\h)$:
\begin{enumerate}
\item $k$ belongs to the norm closure of $\U(h)$;
\item for any open subset $\Os \subset \C$, 
$$\mbox{\textnormal{rank}}( \chi_\Os(h)) = \mbox{\textnormal{rank}}(\chi_\Os(k));$$
\item (when $\h$ is separable) $h$ and $k$ have the same spectrum and the same multiplicity at each of their isolated eigenvalues.
\end{enumerate}
\end{theorem}
Here and throughout, ``rank" means the cardinal dimension of the closure of the range.  The cardinal-valued function defined on the open subsets of $\C$,
$$M_h: \Os \mapsto \text{rank}( \chi_\Os(h)),$$
is called the \textit{crude multiplicity function} for $h$ (\cite{AD}).  (Sometimes this terminology is used for a related function defined on all subsets of $\C$, or on $\C$ itself.  See Remark \ref{T:ptcrude}.)  The crude multiplicity function is therefore a complete invariant for the norm-closed unitary orbit of a normal operator in $\B(\h)$.

Theorem \ref{T:normbh} is conceptually close to Berg's 1971 generalization of the Weyl-von Neumann theorem (see Section \ref{SS:wvnb}).  This was generalized to separable representations of separable noncommutative $C^*$-algebras by Voiculescu (\cite{V1976}) in 1976, and in 1981 Hadwin eliminated the separability hypotheses (\cite{H1981}).  We will say more about this later in the introduction.  Applied to strong* closures of unitary orbits of normal operators, Hadwin's work yields

\begin{theorem} \label{T:strongbh} $($\cite[Proposition 5.3]{H1981}$)$
Let $h$ be a normal operator in $\B(\h)$.  The strong* closure of $\U(h)$ consists of the normal operators $k$ which satisfy 
\begin{equation} \label{E:hadcruder} 
\min \{M_k(\Os), \aleph_0\} \leq \min \{M_h(\Os), \aleph_0\}, \qquad \forall \mbox{\textnormal{ open }} \Os \subset \C.
\end{equation}
\end{theorem}

\noindent The case of separable $\h$, in which the minima above are unnecessary, was treated earlier in \cite[Corollary 3.4 and Theorem 4.3]{H1977}.

\smallskip

The foundation of the present paper is

\begin{theorem} \label{T:main}
Let $h$ and $k$ be operators in a von Neumann algebra $\M$ of arbitrary type and cardinality, with $h$ normal.
\begin{enumerate}
\item The following three conditions are equivalent:
\begin{itemize}
\item $k$ belongs to the strong* closure of $\U(h)$;
\item $k$ is normal and belongs to the strong closure of $\U(h)$;
\item $k$ is normal, and for any open $\Os \subseteq \C$, $\chi_\Os(k)$ is a strong limit of projections Murray-von Neumann equivalent to $\chi_\Os(h)$.
\end{itemize}
\item The following two conditions are equivalent:
\begin{itemize}
\item $k$ belongs to the norm closure of $\U(h)$;
\item $k$ is normal, and for any open $\Os \subseteq \C$, $\chi_\Os(h)$ and $\chi_\Os(k)$ are Murray-von Neumann equivalent.
\end{itemize}
\end{enumerate}
\end{theorem}

\noindent In Section \ref{S:back} we review recent work of the author (\cite{S2005}) which shows how the third condition in (1) extends \eqref{E:hadcruder}. 

Theorem \ref{T:main} represents a relativization of Theorems \ref{T:normbh} and \ref{T:strongbh}.  It is natural to express its content by defining a \textit{crude $\M$-multiplicity function}, and we do this in Section \ref{S:reform}.  Essentially by construction, the crude $\M$-multiplicity function is a complete invariant for the norm-closed unitary orbit in $\M$.  For the strong* closure, domination of crude $\M$-multiplicity functions is the right notion as long as $\M$ is not too big.  Put colloquially, spectral mass can decrease in a strong* limit by being ``pushed out to infinity" -- but at cardinalities higher than $\aleph_0$, spectral mass can increase also.  To handle this we also introduce another $\M$-multiplicity function, along the lines of \eqref{E:hadcruder}, which is just a little cruder.

The moral of all this is not surprising: operator theory inside a von Neumann algebra $\M$ replaces the discrete dimensionality of subspaces by Murray-von Neumann equivalence classes of projections in $\M$.  For $\M$ $\sigma$-finite, Theorem \ref{T:main} can be proved much like its predecessors for $\B(\ell^2)$, but the general case requires a fair bit of extra care.  (See Remark \ref{T:sigf}.)  In fact, the majority of this paper takes its cue from classical operator theory, cooking with three extra ingredients: cardinality, type, and nontrivial center.  The first two of these produce several interesting phenomena, while the third mostly has the effect of making the proofs more complicated.

Let us mention some of the main applications of Theorem \ref{T:main}.  Exactly when the ambient algebra is a direct sum of $\sigma$-finite algebras, the equivalence relation on normal operators, ``equal closed unitary orbit," is independent of the choice of topology from among norm, strong*, and strong (Theorem \ref{T:reform}).  In a properly infinite algebra, the strong closure of the unitary orbit of an arbitrary operator meets the center in the (nonempty) left essential central spectrum of Halpern (Theorem \ref{T:lec}).  This leads to a ``type III Weyl-von Neumann-Berg theorem" (Theorem \ref{T:wvnb'}): for two normal operators in a type III algebra, containment of strong-closed unitary orbits is equivalent to containment of central essential spectra.  And finally, we put considerable effort in Section \ref{S:agree} into characterizing agreements between the various closures of the unitary orbit of a normal operator $h \in \M$.  The results of this are as follows:
\begin{itemize}
\item (Theorem \ref{T:strongnorm}) The norm closure is strong*-closed if and only if either $\M$ is finite or $h$ is central.
\item (Theorem \ref{T:ss}) When $\M$ is a factor with separable predual, the strong* closure is strong-closed if and only if the essential spectrum of $h$ has no interior and does not separate the plane.
\item (Theorem \ref{T:normc}) When $\M$ is a factor, the unitary orbit is norm-closed if and only if $h$ is diagonal with countable essential spectrum, and for each point $\lambda$ in the essential spectrum, the spectral projection for some deleted neighborhood of $\lambda$ is strictly subequivalent to the spectral projection for $\{\lambda\}$.
\end{itemize}

\smallskip

Subcases of Theorem \ref{T:main}(2) which deal with von Neumann algebras other than $\B(\h)$ have appeared in the literature, sometimes implicit in results of a different flavor.  For $\sigma$-finite factors, it follows (with a little work, described in more detail in Section \ref{SS:dist}) from some estimates of Hiai and Nakamura for the distance between the unitary orbits of two normal operators (\cite{HN}).  For self-adjoint operators in a $\text{II}_1$ factor, it can be found in papers by Sunder and Thomsen (\cite{ST}) and Arveson and Kadison (\cite{AK}).  To understand the other sources, we need to rephrase Theorem \ref{T:main} in terms of *-homomorphisms of $C^*$-algebras (always assumed to have a unit).

Let $h$ be a normal operator on the Hilbert space $\h$, and let $\{u_\alpha\}$ be a net of unitaries such that $(\text{Ad }u_\alpha)(h)$ converges strong* to an operator $k$.  Since polynomials in $h, h^*$ are norm-dense in $C^*(h) \subset \B(\h)$, it follows that $\text{Ad }u_\alpha: C^*(h) \to \B(\h)$ converges in the point-strong* topology to a *-homomorphism from $C^*(h)$ onto $C^*(k)$.  So the strong* closure of $\U(h)$ may be identified with the point-strong* closure of the $\B(\h)$-inner *-homomorphisms $C^*(h) \to \B(\h)$.  The norm closure of $\U(h)$ corresponds to the point-norm closure of the $\B(\h)$-inner *-homomorphisms $C^*(h) \to \B(\h)$, which consists of isomorphisms.  The analogous statements for $W^*(h)$ are \textit{not true}.  This is directly reflected in the fact that only open sets play a role in Theorem \ref{T:main}.

It makes no difference to let $h$ be non-normal in the preceding paragraph, or even to replace $C^*(h)$ with an arbitrary $C^*$-algebra.  This is the venue for Voiculescu's theorem (\cite{V1976}) and a wealth of interesting work by Hadwin (including \cite{H1977,H1981,H1987,HL}, the last with Larson).  In these papers the ideal of compact operators plays a key role in the questions and answers, as both compact and noncompact representations enjoy special properties.  See Arveson (\cite{A1977}) for related developments, or Davidson (\cite[Sections II.4-5]{D1996}) for an accessible introduction.  One highlight is the following profound generalization of Theorem \ref{T:normbh}.

\begin{theorem} $($\cite[Theorem 3.14]{H1981}$)$ \label{T:hv}
Let $\pi$ and $\rho$ be unital representations of the $C^*$-algebra $A$ on a fixed Hilbert space $\h$.  Then $\rho$ belongs to the point-norm closure of the unitary orbit of $\pi$ $(=\{(\mbox{\textnormal{Ad }}u) \circ \pi \mid u \in \U(\B(\h))\})$ if and only if
\begin{equation} \label{E:hadrank}
\mbox{\textnormal{rank}}(\pi(a)) = \mbox{\textnormal{rank}}(\rho(a)), \qquad \forall a \in A.
\end{equation}
\end{theorem}

For an even broader context, one may allow ``representations" inside algebras other than $\B(\h)$.  Actually the past ten years have seen a steady stream of research concerning notions of approximate unitary equivalence for *-homomor\-phisms $A \to B$, where both $A$ and $B$ are sufficiently well-behaved $C^*$-algebras (\cite{Da,L1997,L2000,GL,L2002,L2004}).  This industry typically uses $KK$-theory (substituting for notions like ``rank") and has generated some powerful results.  But questions which involve the operator topologies on $\B(\h)$ are off the menu.

In this paper we consider a specific subcase: $A$ is singly-generated and abelian, and $B$ is a von Neumann algebra.  It is the second assumption which imparts a different flavor; the availability of spectral projections and operator topologies keeps us closer to classical operator theory.  As far as we can tell, the intersection of Theorem \ref{T:main}(2) with the  $C^*$-papers mentioned above is limited to the case where $B$ is a finite or type III factor, and even this is somewhat obscure due to the difference in language.  (Should the reader attempt the translation, it might be helpful to keep in mind that $K_1$ vanishes for von Neumann algebras.  It should also be plausible that the second condition of Theorem \ref{T:main}(2) implies an equality for morphisms of $K_0$-groups.)

There is a very recent paper by Ding and Hadwin (\cite{DH}) whose agenda is related to ours.  For them the central issue is the possibility of obtaining a relative version of Theorem \ref{T:hv}, replacing $\B(\h)$ by a von Neumann algebra $B$ and ``equal rank" by Murray-von Neumann equivalence of range projections in $B$.  Given the successful transition from Theorems \ref{T:normbh} and \ref{T:strongbh} to Theorem \ref{T:main}, it may be surprising that such a theorem does \textit{not} hold in general.  Hadwin gave a construction using free entropy in \cite[Corollary 3.5]{H1998}, and here we present an additional

\begin{example}
Let $B$ be a $\text{II}_1$ factor equipped with an automorphism $\rho$ which is not a point-norm limit of inner automorphisms.  For example, we may take the automorphism generated by the flip $(x \otimes y) \mapsto (y \otimes x)$ on $\NN \overline{\otimes} \NN$, where $\NN$ is any non-hyperfinite $\text{II}_1$ factor (\cite[Theorem 5.1]{C}).  Then there is a finite set $\{x_j\} \subset B$ such that no unitary $u$ makes $\max_j \{\|u x_j u^* - \rho(x_j)\|\}$ arbitrarily small.  In other words, $\text{id}$ and $\rho \circ \text{id}$ are not approximately unitarily equivalent *-homomorphisms from $C^*(\{x_j\})$ to $B$.  But $\rho$ must preserve the unique tracial state on $B$, and it is easy to see from this that $\text{id}$ and $\rho \circ \text{id}$ satisfy the ``$B$-version" of \eqref{E:hadrank} (\cite[Lemma 3]{DH}).
\end{example}

Still, Ding and Hadwin proved that a relative Theorem \ref{T:hv} is true whenever $B$ has separable predual and $A$ is a direct limit of direct sums of AF-algebras tensored with commutative algebras (\cite[Corollary 3]{DH}).  This result of course encompasses Theorem \ref{T:main}(2) for von Neumann algebras with separable predual (a condition which, as we have noted, allows for some simplification in the proof).  The reader is referred to \cite{H1998,DH} for interesting complementary results; the only other overlap between these papers and the present one is in our Theorem \ref{T:reform}(4), which is more or less contained in \cite[Corollary 11]{DH}.

In fact we have not imitated the \textit{methods} of Voiculescu, Hadwin, or Arveson, and the underlying reason is the unavailability of compact operators.  (But there are related questions, and answers, which do make use of ideals in von Neumann algebras.)  Instead we benefit from the particularly simple form of our initial algebra.  It allows us effectively to suppress representation theory and state the main results in terms of a familiar and concrete object, the spectral measure of a normal operator.

\smallskip

A few words are in order about the choice of topologies.  Our use of strong, strong*, and weak operator topologies might suggest that $\M$ must be given as a subalgebra of some $\B(\h)$.  This is not the case.  These topologies do depend on the choice of representation, but not when restricted to bounded subsets of $\M$, and unitary orbits are obviously bounded.  The reader who finds the semantics awkward can replace the three topologies with the $\sigma$-strong, $\sigma$-strong*, and $\sigma$-weak; these last admit a nonspatial definition and agree with their counterparts on bounded sets.   

Except in parts of Sections \ref{S:esssp} and \ref{S:agree}, the operators in this paper are normal and described in terms of spectral measures.  The strong* topology is therefore suitable because the strong* limit of normal operators is again normal.  The same is not true of the strong and weak topologies.  (The strong limits of normal operators are exactly the subnormal operators (\cite{Bi,CoH}), which play a role in Section \ref{S:agree}.  In the context of *-homomorphisms of $C^*$-algebras discussed above, the point-strong and point-weak topologies turn out to be equal to the point-strong* anyway (\cite[Section II.4]{D1996}).)  Yet the strong* and strong topologies do agree when \textit{restricted} to the set of normal operators (\cite[Proposition II.4.1]{T}); this means that the strong*-closed unitary orbit of a normal operator is precisely the set of normal operators in the strong-closed unitary orbit.  In our arguments we typically use the simpler strong topology when we have assumed normality of the limit.

\smallskip

The paper is organized as follows.  The next section reviews all the ingredients for the proof of Theorem \ref{T:main}, which is accomplished in Section \ref{S:proof}.  Then we give a short section of corollaries.  In Section \ref{S:reform} we reformulate results in terms of crude and cruder $\M$-multiplicity functions (Theorem \ref{T:reform}).  Section \ref{S:esssp} discusses various essential spectra and their relation to central elements in closed unitary orbits of arbitrary operators.  Then we spend a section briefly discussing the relationship between our results and two topics with substantial literature: Weyl-von Neumann-Berg theorems and distance between unitary orbits of normal operators.  Finally, Section \ref{S:agree} is the longest of the paper.  It considers the obstacles to a description of $\U(h)$ itself, and derives necessary and sufficient conditions for agreements between $\U(h)$ and its closures.

There are many possible directions for continuing this line of research.  At the least there will be a sequel paper, joint with C. Akemann, describing \textit{weak} closures of unitary orbits.  These sets are often convex - a fact which seems to have been overlooked - and are quite different from the closed orbits of this paper.

\section{Background} \label{S:back}
Let $\M$ be a von Neumann algebra of arbitrary type and cardinality.  Along with the usual notations we write $\z(\M)$ for its center and $\p(\M)$ for its projection lattice.  Where appropriate, we symbolize the norm, strong*, strong, and weak operator topologies by $\|\|$, $s^*$, $s$, and $w$.  The support of a self-adjoint operator is $s(\cdot)$, while the central support and spectrum of an arbitrary operator are denoted $c(\cdot)$ and $\text{sp}(\cdot)$, respectively.  The open disk of radius $r$ centered at $\lambda \in \C$ is written $B_r(\lambda)$.  For a set $A \subseteq \C$, $A^c$ is the complement of $A$ and $\partial A$ its boundary.  Finally, a normal operator $h$ in a von Neumann algebra is said to be \textit{diagonal} if its spectral measure is atomic, i.e. $\sum_{\lambda \in \C} \chi_{\{\lambda\}}(h) = 1$.

We use the standard terminology and results from \cite[Section V.1]{T} for projections, including $p^\perp$ for $1-p$.  Our notations are $p \sim q$ for (Murray-von Neumann) equivalence, $p \preccurlyeq q$ for subequivalence, and $p \prec q$ for $p \preccurlyeq q$ but not $p \sim q$. 
Note that for pairwise orthogonal sets $\{p_\alpha\}$, $\{q_\alpha\}$,
\begin{equation} \label{E:addeq}
p_\alpha \sim q_\alpha, \: \forall \alpha \quad \Rightarrow \quad \left(\sum p_\alpha\right) \sim \left(\sum q_\alpha\right),
\end{equation}
\begin{equation} \label{E:addeq2}
p_\alpha \preccurlyeq q_\alpha, \: \forall \alpha \quad \Rightarrow \quad \left(\sum p_\alpha\right) \preccurlyeq \left(\sum q_\alpha\right).
\end{equation}

One term which may cause confusion is \textit{properly infinite}; it describes a nonzero projection $p$ for which $zp$ is infinite or zero for any central projection $z$.  Any adjective can be applied to an algebra when the adjective describes the identity projection of the algebra.

A properly infinite projection $p$ can be decomposed into a (finite or infinite) countable sum of projections, each of which is equivalent to $p$.  We record two simple consequences.

\begin{lemma} \label{T:add}
Let $\{p,q,q_j\} \subset \p(\M)$, with $\{q_j\}$ countable and $q$ properly infinite.
\begin{enumerate}
\item If $p \preccurlyeq q$ and $p \perp q$, then $q \sim (p+q)$.
\item If $q \sim q_j$ for all $j$, then $q \sim \vee q_j$.
\end{enumerate}
\end{lemma}

\begin{proof}
For the first part, write $q = q_1 + q_2$, where $q_1 \sim q_2 \sim q$.  By \eqref{E:addeq2},
$$ q \preccurlyeq (p + q) \preccurlyeq (q_1 + q_2) = q.$$
The second part (also easy) is Lemma 3.2(1) in \cite{S2005}.
\end{proof} 

By \eqref{E:addeq}, one may unambiguously sum any set in $\ps$ for which there are mutually orthogonal representatives, simply by taking the equivalence class of the sum of representatives.  This endows $\ps$ with a partially-defined addition making the quotient map $\p(\M) \twoheadrightarrow \ps$ completely additive (\cite[Section 5]{S2005}).  It also determines an order on $(\p(\M)/\sim)$: $[p] \leq [q]$ if there exists a projection $r$ with $[p] + [r] = [q]$.  One may induce the same order directly from the quotient operation, i.e.
$$\exists q_1, q_2 \text{ with } q_1 \sim p_1, \: q_2 \sim p_2, \: q_1 \leq q_2 \iff [p_1] \leq [p_2].$$
So $[p_1] \leq [p_2]$ means nothing other than $p_1 \preccurlyeq p_2$, and $(\p(\M)/\sim)$ is totally ordered if and only if $\M$ is a factor.

We showed in \cite[Theorem 6.1]{S2005} that $\ps$ is a complete lattice.  The proof of this is technical, but lattice operations on pairs are easily expressed via the comparison theorem for projections (\cite[Theorem V.1.8]{T}).  For $p,q \in \p(\M)$, let $z$ be a central projection with $zp \preccurlyeq zq$, $z^\perp p \succcurlyeq z^\perp q$.  Then 
\begin{equation} \label{E:lattice}
[p] \wedge [q] = [zp + z^\perp q], \qquad [p] \vee [q] = [z^\perp p + zq].
\end{equation}
It follows from \eqref{E:lattice} and Lemma \ref{T:add}(1) that when $\{p_j\}_{j=1}^n \subset \p(\M)$ and $\sum p_j$ is properly infinite,
\begin{equation} \label{E:easyvee}
\left[\sum p_j \right] = \vee [p_j].
\end{equation}

In line with the terminology for linear maps, we say that the quotient map $\p(\M) \twoheadrightarrow \ps$ is \textit{normal} if $[\vee p_\alpha] = \vee [p_\alpha]$ for all increasing nets $\{p_\alpha\} \subset \p(\M)$.

\begin{proposition} \label{T:normal} $($\cite[Proposition 5.1]{S2005}$)$
The quotient map $\p(\M) \twoheadrightarrow \ps$ is normal if and only if $\M$ is a (possibly uncountable) direct sum of $\sigma$-finite algebras.
\end{proposition}

The reverse implication of Proposition \ref{T:normal} is explained by the presence of an extended center-valued trace which faithfully parameterizes $\ps$.  This is a map $T$ from the positive cone of a von Neumann algebra to the extended positive cone of its center; see \cite[Chapter V.2]{T} for details.  

Semifinite von Neumann algebras are characterized by the existence of a faithful normal semifinite extended center-valued trace $T$ (\cite[Theorem V.2.34]{T}).  If $\M$ is finite, there is a \textit{unique} faithful extended center-valued trace $T$ with $T(1) = 1$ (\cite[Theorem V.2.6]{T}).  Such a map is automatically normal and $\sigma$-strong-$\sigma$-strong continuous (\cite[Theorem 13]{G}), and the linear extension which is defined on all of $\M$ is called simply a \textit{center-valued trace}.

\begin{convention}
Whenever we talk of an  ``extended center-valued trace" $T$ on $\M_+$ in the sequel, it is assumed that
\begin{itemize}
\item $T$ is normal and faithful;
\item on the finite summand, $T$ agrees with the center-valued trace;
\item on the semifinite summand, $T$ is semifinite;
\item on the infinite type I summand, $T$ takes an abelian projection to its central support.  
\end{itemize}
So $T(h) = (+\infty) c(h)$ if $h \in \M_+$ is supported on the type III summand.
\end{convention}

It can be shown that the normality of $T$ implies weak lower-semicont\-inuity, but to make this precise requires discussion of the complete lattice structure of the extended positive cone of the center.  The following fact will be a sufficient substitute for our purposes.

\begin{lemma} \label{T:liminf} $($\cite[Lemma 3.1]{S2005}$)$
Let $\{x_\alpha\}$ be a net in a semifinite von Neumann algebra $\M$ equipped with an extended center-valued trace $T$.  If $x_\alpha^* x_\alpha = y_1$ is fixed, while $x_\alpha x_\alpha^* \overset{w}{\to} y_2$, then $T(y_1) \geq T(y_2)$.
\end{lemma}

The quotient $\ps$ is sometimes called the \textit{dimension theory} for $\M$.  Extended center-valued traces provide an avenue for parameterizing $\ps$, and in this capacity they are called \textit{dimension functions}.  Exactly when $\M$ is a direct sum of $\sigma$-finite algebras, we have
\begin{equation} \label{E:par}
p \preccurlyeq q \iff T(p) \leq T(q), \qquad p,q \in \p(\M).
\end{equation}
See \cite[Proposition V.1.39 and Corollary V.2.8]{T} and \cite[Proposition 3.7]{S2005}.  (In \cite{S2005} we construct \textit{fully extended center-valued traces}, for which \eqref{E:par} is valid with no restriction on $\M$.)  Keeping \eqref{E:par} in mind, here are some simple examples of dimension theory.

\begin{enumerate}
\item When $\M$ is a type $\text{I}_\kappa$ factor, $\ps$ is isomorphic to the initial segment of cardinals $\leq \kappa$, via the map that sends a projection to its rank.
\item When $\M$ is a type $\text{II}_1$ factor, $\ps \simeq [0,1]$.
\item When $\M$ is a $\sigma$-finite type $\text{II}_\infty$ factor, $\ps \simeq [0, +\infty]$.
\item When $\M$ is a $\sigma$-finite type $\text{III}$ factor, $\ps \simeq \{0, +\infty\}$.
\end{enumerate}

We endow $\ps$ with the quotient of the strong (equivalently, the weak) topology on $\p(\M)$.  Abbreviating this ``quotient operator topology" to ``$QOT$," we have the following description for closures of singletons.

\begin{theorem} \label{T:singleton} $($\cite[Theorems 3.3 and 3.4]{S2005}$)$
Let $p$ be a projection in a von Neumann algebra $\M$.
\begin{enumerate}
\item If $\M$ is finite, the center-valued trace $T$ implements a homeomorphism between $(\ps, QOT)$ and a topological subspace of $(\z(\M)_1^+, \text{\textnormal{strong}})$.  Consequently
\begin{equation} \label{E:finclos}
\overline{\{[p]\}}^{QOT} = \{[p]\}.
\end{equation}
\item If $\M$ is properly infinite and $p$ is finite,
\begin{equation} \label{E:infclos1}
\overline{\{[p]\}}^{QOT} = \{[q] \mid [q] \leq [p]\}.
\end{equation}
\item If $p$ is properly infinite and $c(p)=I$,
\begin{equation} \label{E:infclos2}
\overline{\{[p]\}}^{QOT} = (\p(\M)/\sim).
\end{equation}
\end{enumerate}
Equations \eqref{E:infclos1} and \eqref{E:infclos2} can be synthesized as follows: for $\M$ properly infinite with $T$ an extended center-valued trace,
\begin{equation} \label{E:infclos3}
\overline{\{[p]\}}^{QOT} = \{[q] \mid T(q) \leq T(p) \}, \qquad p \in \p(\M).
\end{equation}
\end{theorem}

We abbreviate the condition $[q] \in \overline{\{[p]\}}^{QOT}$ to $q \in p$.  It simply means that $q$ is a strong limit of projections equivalent to $p$, and will be heavily used in the sequel.  Viewing $\M = \B(\ell^2)$ as the base case, we think of ``$\in$" as a generalization of having equal or lower rank.  For general infinite-dimensional $\B(\h)$, ``$\in$" is the relation expressed in \eqref{E:hadcruder}.

\section{Proof of Theorem \ref{T:main}} \label{S:proof}

We mentioned in the introduction that the strong and strong* topology agree on the set of normal operators.  To establish Theorem \ref{T:main}, it therefore suffices to to assume that $h$ and $k$ are normal operators in a von Neumann algebra $\M$, and prove
\begin{equation} \label{E:main}
k \in \overline{\U(h)}^s \iff \chi_\Os(k) \in \chi_\Os(h), \: \forall \text{ open } \Os \subseteq \C,
\end{equation}
\begin{equation} \label{E:main2}
k \in \overline{\U(h)}^{\|\|} \iff \chi_\Os(k) \sim \chi_\Os(h), \: \forall \text{ open } \Os \subseteq \C.
\end{equation}

Since all four conditions on $h$ and $k$ are compatible with central decompositions, when convenient we may work inside the different type summands of $\M$.  We consider each of the implications in \eqref{E:main} and \eqref{E:main2} separately.

\smallskip

\textit{Proof of \eqref{E:main}, $(\Rightarrow)$:} By Theorem \ref{T:singleton} and restriction to a central summand, it suffices to show that $\chi_\Os(k) \in \chi_\Os(h)$ when $\chi_\Os(h)$ is finite with full central support.  In this case $\M$ is semifinite; let $T$ be an extended center-valued trace on $\M$.  

Define the sequence of continuous functions
$$f_n(\cdot) = \min\{1, n \: \text{dist}(\cdot, \Os^c)\}.$$
These are ``plateau functions" which increase pointwise to $\chi_\Os$.  Each $f_n$ is strong-strong continuous as an operator function by \cite[Theorem II.4.7]{T}.  Therefore $u_\alpha h u_\alpha^* \overset{s}{\to} k$ implies
\begin{equation} \label{E:plateau}
u_\alpha f_n(h) u_\alpha^* = f_n(u_\alpha h u_\alpha^*) \overset{s}{\to} f_n(k).
\end{equation}

We now apply Lemma \ref{T:liminf}, with $x_\alpha = u_\alpha f_n(h)^{1/2}$, concluding
\begin{equation} \label{E:trace}
T(f_n(k)) \leq T(f_n(h)).
\end{equation}
As $n$ increases, the normality of $T$ provides
\begin{equation} \label{E:trace2}
T(\chi_\Os(k)) \leq T(\chi_\Os(h)),
\end{equation}
and by finiteness of $\chi_\Os(h)$ we get $\chi_\Os(k) \preccurlyeq \chi_\Os(h)$.  On the finite summand of $\M$, the $\sigma$-strong continuity of $T$ means that \eqref{E:trace} and \eqref{E:trace2} are equalities, so actually $\chi_\Os(k) \sim \chi_\Os(h)$.  We conclude from Theorem \ref{T:singleton} that $\chi_\Os(k) \in \chi_\Os(h)$.

\smallskip

\textit{Proof of \eqref{E:main2}, $(\Rightarrow)$:} The relation on the left hand side of $\eqref{E:main2}$ is an equivalence relation.  In particular it is symmetric and stronger than $k \in \overline{\U(h)}^s$, so by the previous argument, $\chi_\Os(k) \in \chi_\Os(h)$ and $\chi_\Os(h) \in \chi_\Os(k)$.  A look at Theorem \ref{T:singleton} shows that $\chi_\Os(h) \sim \chi_\Os(k)$ whenever either is finite.

So assume $\chi_\Os(h)$ properly infinite, and let $\{u_j\}$ be a sequence of unitaries such that $u_j h u_j^* \overset{\|\|}{\to} k$.  Set $f(\cdot) = \text{dist}(\cdot, \Os^c)$.  As before, we have that
\begin{equation} \label{E:sconv}
u_j f(h) u_j^* = f(u_j h u_j^*) \overset{s}{\to} f(k).
\end{equation}
(The convergence is actually in norm (\cite[Proposition I.4.10]{T}), but we do not need this.) 

Now we use an obvious fact about supports and Lemma \ref{T:add}(2) to compute
$$\chi_\Os(k) = s(f(k)) \leq \bigvee_j s(u_j f(h) u_j^*) = \bigvee_j u_j \chi_\Os(h) u_j^* \sim \chi_\Os(h).$$
Since the roles of $h$ and $k$ can be reversed, we are done.

\smallskip

\textit{Proof of \eqref{E:main}, $(\Leftarrow)$:} We first consider the special case where both $h$ and $k$ are simple operators satisfying the right-hand side of \eqref{E:main}.  This means that the spectral measures are supported on a finite set of points $\{\lambda_i\}$.  We abbreviate $\chi_{\{\lambda_i\}}(\cdot)$ to $\chi_i(\cdot)$.  When $\M$ is finite this immediately implies that $k = uhu^*$ for some unitary $u$, so we take up the properly infinite case.

Since the identity is infinite and equal to the finite sum $\sum \chi_i(k)$, by \eqref{E:easyvee} we can find central projections $z_i$ with sum 1 and $z_i \chi_i(k) \sim z_i$.  Similarly, we can find central projections $y_i$ with sum 1 and $y_i \chi_i(h) \sim y_i$.  It suffices to explain the construction for each $(z_{i_1} \wedge y_{i_2}) \M$, since they can be assembled into a global solution.

Let us ease the notation by writing $(z_{i_1} \wedge y_{i_2}) \M$ as $\M$, so that $\chi_{i_1}(k) \sim \chi_{i_2}(h) \sim 1$.  The goal is to find a net of unitaries $\{u_\alpha\}$ so that $u_\alpha \chi_i(h) u_\alpha^* \overset{s}{\to} \chi_i(k)$ for each index $i$.  By taking a finite linear combination, this implies $u_\alpha h u_\alpha^* \overset{s}{\to} k$.

For each $i \neq i_1, i_2$, we use the hypothesis $\chi_i(k) \in \chi_i(h)$ to find a subprojection $p_i \leq \chi_i(h)$ with the following properties.  Let $\chi_i(k) = z \chi_i(k) + z^\perp \chi_i(k)$ be a decomposition into finite and properly infinite parts.  We require $zp_i \sim z\chi_i(k)$, $z^\perp p_i$ properly infinite, $c(z^\perp p_i) = c(z^\perp \chi_i(k))$, and $z^\perp p_i \preccurlyeq z^\perp \chi_i(k)$.  By Theorem \ref{T:singleton}(3) we can find a net of partial isometries $\{v_{\alpha_i}\}$ between $p_i$ and subprojections of $\chi_i(k)$ which increase to $\chi_i(k)$.  (Here $\alpha_i$ belongs to an index set which depends on $i$.)  For $i=i_1$, we do the same thing with the additional requirement that $[\chi_{i_1}(k) - v_{\alpha_{i_1}} v_{\alpha_{i_1}}^*] \sim 1$ - recall that $\chi_{i_1}(k) \sim 1$.  For $i=i_2$, we do the same thing with the additional requirement that $(\chi_{i_2}(h) - p_{i_2}) \sim 1$ - recall that $\chi_{i_2}(h) \sim 1$.  

We consider the product net of the nets so far defined.  An element is a choice of $\alpha_i$ for each $i$, corresponding to the partial isometry
$\sum_i v_{\alpha_i}$ between $\sum_i p_i$ and $\sum_i v_{\alpha_i} v_{\alpha_i}^*$.  The complementary projections are
\begin{equation} \label{E:lr}
\sum_i (\chi_i(h) - p_i) \quad \text{and} \quad \sum_i (\chi_i(k) - v_{\alpha_i}v_{\alpha_i}^*).
\end{equation}
We know that $(\chi_{i_2}(h) - p_{i_2})$ and $(\chi_{i_1}(k) - v_{\alpha_{i_1}} v_{\alpha_{i_1}}^*)$ are both $\sim 1$, so the two projections of \eqref{E:lr} are equivalent.  Choose any partial isometry between them, and add it to $\sum_i v_{\alpha_i}$ to get a unitary $u_{\alpha_i}$.  In this way we have constructed a net of unitaries $\{u_{\alpha_i}\}$, indexed by the choices of $(\alpha_i)$ and ordered by the product order.  By construction we have the strong limits
\begin{align*}
\lim u_{\alpha_i}\chi_i(h)u_{\alpha_i}^* &= \lim u_{\alpha_i} p_i u_{\alpha_i}^* + \lim u_{\alpha_i}(\chi_i(h) - p_i) u_{\alpha_i}^* \\ &= \lim v_{\alpha_i} p_i v_{\alpha_i}^* + \lim u_{\alpha_i}(\chi_i(h) - p_i) u_{\alpha_i}^* \\ &= \chi_i(k) + \lim u_{\alpha_i}(\chi_i(h) - p_i) u_{\alpha_i}^*.
\end{align*}
Now notice that from the equivalence in \eqref{E:lr},
$$ u_{\alpha_i}(\chi_i(h) - p_i) u_{\alpha_i}^* \leq \sum (\chi_i(k) - v_{\alpha_i} v_{\alpha_i}^*) \overset{s}{\to} 0.$$
Therefore we have shown that $u_{\alpha_i}\chi_i(h) u_{\alpha_i}^* \overset{s}{\to} \chi_i(k)$ for all $i$.  This completes the argument for $(z_{i_1} \wedge y_{i_2}) \M$, and by gluing finishes the proof of the finite spectrum case.

Now take arbitrary normal $h,k$ satisfing the right-hand side of \eqref{E:main}.  We construct a sequence of rectangular grids $\{\G^j\}$ having certain properties:
\begin{itemize}
\item each $\G^j \subset \C$ is a closed rectangle with vertical and horizontal sides, gridded into finitely many smaller rectangles;
\item $\text{sp}(h)$ (which contains $\text{sp}(k)$, by taking $\Os = \text{sp}(h)^c$ in \eqref{E:main}) is in the interior of each $\G^j$;
\item the mesh of $\G^j$ (the largest diameter of a subrectangle) goes to 0 as $j \to \infty$;
\item all vertical line segments in the $\G^j$ are disjoint, as are all horizontal line segments.
\end{itemize}
We consider the grid $\G^j$ to be the disjoint union of open rectangles $\{R^j\}$, the open line segments $\{S^j$\} which form the sides of the rectangles, and the vertices $\{V^j\}$ of the grid.  It would be proper to index elements of these sets with subscripts, but given the preponderance of indices, we will do this only when necessary.  The idea is to use the grids to produce approximants to $h$ and $k$ with finite spectrum.

First, on the finite summand of $\M$, the assumption implies that $\chi_{R^j}(k) \sim \chi_{R^j}(h)$ for each open rectangle $R^j$.  For any $S^j$, we may consider the open rectangle formed by $S^j$ and the two adjacent rectangles $R^j_1, R^j_2$:
$$(\chi_{R^j_1}(k) + \chi_{S^j}(k) + \chi_{R^j_2}(k)) \sim (\chi_{R^j_1}(h) + \chi_{S^j}(h) + \chi_{R^j_2}(h)).$$
From the previous statement and the finiteness of $\M$, it follows that $\chi_{S^j}(k) \sim \chi_{S^j}(h)$.  Similarly, for any $V^j$, we consider the open rectangle formed by $V^j$, four segments, and four rectangles: from the previous statements it follows that $\chi_{V^j}(h) \sim \chi_{V^j}(k)$.

Now for each $j$, let $\{\lambda^j\}$ be the (finitely many) center points of the $R^j$.  For each of the $\lambda^j$, let $T^j$ be the union of the rectangle containing $\lambda^j$, the right and bottom sides of the rectangle, and the right bottom corner.  Consider the approximations
$$h^j = \sum \lambda^j \chi_{T^j}(h), \qquad k^j = \sum \lambda^j \chi_{T^j}(k),$$
where the sums are taken over the center points $\lambda^j$.  Then $h^j$ and $k^j$ have finite spectrum and satisfy the right-hand side of \eqref{E:main}; in particular there is a unitary $u^j$ with $k^j = u^j h^j {u^j}^*$.  By the mesh condition, $\|h^j - h\|, \|k^j - k\| \to 0$.  So
$$\|u^j h {u^j}^* - k\| \leq \|u^j h {u^j}^* - u^j h^j {u^j}^* \| + \| k^j - k \| \overset{j \to \infty}{\to} 0.$$
This shows that $k \in \overline{\U(h)}^{\|\|} \subseteq \overline{\U(h)}^s$.

We now assume that $\M$ is properly infinite.  We form $h^j$ exactly as before, but we will need a little more care with $k^j$ because the spectral projections for $h$ and $k$ corresponding to segments and vertices do not follow a simple relation.

Since $1 = \sum \chi_{\lambda^j}(h^j)$ (sum over $\lambda^j$), it follows from \eqref{E:easyvee} that one can find finitely many central projections $\{z^j_n\}$, with $\sum_n z^j_n = 1$, such that each $z_n$ is equivalent to some $\chi_{\lambda_n^j}(h^j)$.  We form the approximation
$$k^j = \left(\sum_{R^j} \lambda^j \chi_{R^j}(k) \right) + \chi_{\cup S^j \cup V^j}(k) \left( \sum_n \lambda^j_n z_n \right).$$
In other words, we approximate the open rectangles by their center points, and everything else we move, summand by summand, to a center point where $h^j$ is properly infinite.  It follows that $k^j$ and $h^j$ satisfy the right-hand side of \eqref{E:main}.  Letting $q^j$ be the projection $\chi_{\cup S^j \cup V^j}(k)$, we have that $k^j$ is a small norm perturbation of $k$ off $q^j$.

We also have that $q^j$ goes to 0 strongly.  To see this, notice that while the segments and vertices of any two $\G^j$ may intersect nontrivially, our assumptions imply that there are no points which belong to the segments and vertices of three different $\G^j$.  Therefore the $q^j$ are mutually commuting projections such that the product of any three is zero.  We set
$$q^{ij} = q^i q^j, \quad i < j; \qquad q^{ii} = q^i - \left(\sum_{i<j} q^{ij}\right) -  \left(\sum_{k<i} q^{ki}\right).$$
Then $\{q^{ij}\}_{i \leq j}$ are pairwise disjoint, and
$$q^n \leq \sum_{j \geq n} q^{ij} \overset{s}{\to} 0.$$

For any $u \in \U(\M)$, we have
$$uhu^* - k = (uhu^* - uh^j u^*) + (uh^ju^* - k^j) + (k^j - k)(1-q^j) + (k^j - k)q^j.$$
In the sum, the first and third terms have norm $\leq \text{mesh} (\G^j)/2$.  Choosing large enough $j$, then an appropriate $u$ from the first part of the proof, will put the second and fourth terms in any strong neighborhood of 0.  It follows that $k \in \overline{\U(h)}^s$.

\smallskip

\textit{Proof of \eqref{E:main2}, $(\Leftarrow)$:} We have already shown this for the finite summand of $\M$, so we assume that $\M$ is properly infinite.  Let $\{\G^j\}$, $\{R^j\}$, $\{\lambda^j\}$, etc., be as before.  Instead of directly defining approximants for $h$ and $k$ based on $\G^j$, for each $j$ we explain a finite series of norm-small changes which ends up at approximants whose spectra are contained in $\{\lambda^j\}$.  The obstacle is again that spectral projections for $h$ and $k$ corresponding to an $S^j$ or $V^j$ need not be equivalent.

The index $j$ is fixed in the next three paragraphs; we continue to omit additional subscripts unless necessary.  Note that (*) \textit{each change preserves $\chi_\Os(k) \sim \chi_\Os(h)$ for all open sets which are unions of the $\{R^j\}$, $\{S^j\}$, and $\{V^j\}$}.

First replace the spectral measure for each rectangle $R^j$ by a point mass at the center $\lambda^j$.  In other words, add $\sum_{R^j} (\lambda^j - h) \chi_{R^j}(h)$ to $h$, and similarly for $k$.  We understand this kind of construction whenever we speak of ``moving" a projection or set below.

For each segment $S^j_0$ which is interior in $\G^j$, we apply the right-hand side of \eqref{E:main2} to the open rectangle formed by $S^j_0$ and the two adjacent rectangles $R^j_1, R^j_2$.  We have
\begin{equation} \label{E:someequiv1}
\chi_{R^j_1}(h) \sim \chi_{R^j_1}(k), \quad \chi_{R^j_2}(h) \sim \chi_{R^j_2}(k),
\end{equation}
\begin{equation} \label{E:someequiv2}
(\chi_{S^j_0}(h) + \chi_{R^j_1}(h) + \chi_{R^j_2}(h)) \sim (\chi_{S^j_0}(k) + \chi_{R^j_1}(k) + \chi_{R^j_2}(k)).
\end{equation}
On the central summand where the projections of \eqref{E:someequiv2} are finite, necessarily $\chi_{S^j_0}(h) \sim \chi_{S^j_0}(k)$, and we can move both to either adjacent $\lambda^j$.  On the complementary summand, the projections of \eqref{E:someequiv2} are properly infinite and \eqref{E:easyvee} applies.  So we may divide this into three central subsummands, on each of which we have one of the following:
\begin{enumerate}
\item $(\chi_{S^j_0}(h) + \chi_{R^j_1}(h) + \chi_{R^j_2}(h)) \sim \chi_{R^j_1}(h)$;
\item $(\chi_{S^j_0}(h) + \chi_{R^j_1}(h) + \chi_{R^j_2}(h)) \sim \chi_{R^j_2}(h)$;
\item $(\chi_{S^j_0}(h) + \chi_{R^j_1}(h) + \chi_{R^j_2}(h)) \sim \chi_{S^j_0}(h)$.
\end{enumerate}
We assume that the central summands for (1) and (2) have been taken as large (together) as possible.  Then on summand (3), $\chi_{R^j_1}(h)$ and $\chi_{R^j_2}(h)$ are strictly subequivalent to $\chi_{S^j_0}(h) + \chi_{R^j_1}(h) + \chi_{R^j_2}(h)$.  This implies that the spectral projections also satisfy the same conditions for $k$, and in particular $\chi_{S^j_0}(h) \sim \chi_{S^j_0}(k)$ on (3).  Now on (1) and (3), move $S^j_0$ to $\lambda^j_1$; on (2), move $S^j_0$ to $\lambda^j_2$.  By Lemma \ref{T:add}(1), condition (*) is maintained, no matter the order in which we move the $S^j$.

Now for each interior $V^j_0$, we consider the open rectangle composed of $V^j_0$, the four adjacent $S^j$, and the four adjacent $R^j$.  Spectral projections for the $S^j$ have become zero.  By a procedure analogous to the preceding paragraph, we move $\chi_{V^j_0}(\cdot)$ to the adjacent $\lambda^j$, preserving (*).

This process ends with norm approximants $h^j$ and $k^j$, each of which has spectrum in $\{\lambda^j\}$.  Since $\chi_{\lambda^j}(k^j) \sim \chi_{\lambda^j}(h^j)$ for each $\lambda^j$, there is a unitary $u^j$ with $u^j h^j {u^j}^* = k^j$.  Finally,
$$\|u^j h {u^j}^* - k\| \leq \|u^j h {u^j}^* - u^j h^j {u^j}^*\| + \|k^j - k \| \leq \text{mesh}(\G^j) \to 0.$$
This proves that $k \in \overline{\U(h)}^{\|\|}$.

\begin{remark} \label{T:sigf}
If $\M$ is $\sigma$-finite, the proof of the last two implications above can be simplified significantly by choosing $\{\G^j\}$ so that spectral projections for all $S^j$ and $V^j$ are zero.  Gellar-Page (\cite{GP}) and Ding-Hadwin (\cite{DH}) take advantage of this.
\end{remark}

\begin{remark} \label{T:noclosed}
A version of Theorem \ref{T:main} involving non-open sets will not work.  Set $h$ to be the diagonal operator in $\B(\ell^2)$ with diagonal $(1, 1/2, 1/3, \dots)$.  Let $u_j$ be the permutation matrix corresponding to the cycle $(12 \cdots j)$.  Then $u_j h u_j^*$ is the diagonal $(1/j, 1, 1/2, \dots 1/(j-1), 1/(j+1), \dots)$ and converges in norm to the diagonal $(0, 1, 1/2, \dots) = k$.  Thus $\chi_{\{0\}}(k)$ has rank one, while $\chi_{\{0\}}(h)$ is zero.
\end{remark}

\section{Some corollaries of Theorem \ref{T:main}} \label{S:cor}

Although we use Theorem \ref{T:main} to prove the statements in this section, some may also be shown directly.  It will be useful to make them explicit.
  
\begin{corollary} \label{T:spinc}
Let $h,k \in \M$ be normal.
\begin{equation} \label{E:spinc1}
k \in \overline{\U(h)}^{\|\|} \Rightarrow \mbox{\textnormal{sp}}(k) = \mbox{\textnormal{sp}}(h).
\end{equation}
\begin{equation} \label{E:spinc2}
k \in \overline{\U(h)}^s \Rightarrow \mbox{\textnormal{sp}}(k) \subseteq \mbox{\textnormal{sp}}(h).
\end{equation}
\end{corollary}

\begin{proof} The implication \eqref{E:spinc2} was mentioned in the last section.  (Set $\Os = \text{sp}(h)^c$ in \eqref{E:main}.)  When $k$ belongs to $\overline{\U(h)}^{\|\|} \subseteq \overline{\U(h)}^s$, then $h$ belongs to $\overline{\U(k)}^{\|\|} \subseteq \overline{\U(k)}^s$ as well, so the two spectra must be equal.
\end{proof}

\begin{corollary} \label{T:converse}
Let $h,k \in \M$ be normal.
\begin{enumerate}

\item Assume either that $\M$ is a $\sigma$-finite type III factor, or that $\M$ is a $\sigma$-finite type I factor and $\mbox{\textnormal{sp}}(h)$ contains no isolated points.  Then the converse to \eqref{E:spinc1} holds.

\item Assume either that $\M$ is a type III factor, or that $\M$ is a type I factor and $\mbox{\textnormal{sp}}(h)$ contains no isolated points.  Then the converse to \eqref{E:spinc2} holds.

\end{enumerate}
\end{corollary}
\begin{proof}
Under the hypotheses of (1), any nonzero spectral projection for an open set is equivalent to the identity.  Under the hypotheses of (2), any nonzero spectral projection for an open set is infinite.  Apply \eqref{E:main} and \eqref{E:main2}.
\end{proof}
Theorem \ref{T:wvnb'} gives the non-factor version of the type III part of statement (2).

\smallskip

In the rest of this section we compare $(\p(\M)/\sim)$ with $(\p(\M)/\sim_u)$.  Here $\sim_u$ denotes unitary equivalence: $p \sim_u q$ when there is a unitary $u$ with $upu^* = q$.  To start with, $(\p(\M)/\sim_u)$ does inherit a partial order from $\p(\M)$.  The only nontrivial point is antisymmetry: if $p,q$ are projections and $u,v$ are unitaries with $upu^* \leq q$ and $vqv^* \leq p$, then $p \preccurlyeq q \preccurlyeq p$ and so $p \sim q$.  But also $(1-p) \preccurlyeq (1-q) \preccurlyeq (1-p)$, so $(1-p) \sim (1-q)$.  Together these imply $p \sim_u q$.  As an example, take $\M = \B(\ell^2)$: unitary equivalence classes of projections are identified by giving rank and corank, and the order looks like $(0, \infty), (1,\infty), \dots (\infty, \infty), \dots (\infty, 1), (\infty, 0).$  In fact it is easy to check that \eqref{E:lattice} gives the lattice operations in $(\p(\M)/\sim)$ as well.

However, the analogues of \eqref{E:addeq} and \eqref{E:addeq2} are \textit{not} generally true for unitary equivalence, so $(\p(\M)/\sim_u)$ does not carry a well-defined sum operation.  In $(\B(\ell^2)/\sim_u)$, the sum $(\infty, \infty) + (\infty, \infty)$ could be any $(\infty, c)$.  (Note that Murray-von Neumann equivalence implies unitary equivalence in $M_2 \bar{\otimes} \M$, so the two are more or less interchangeable in $K$-theory (\cite{W-O}).)  Here is the analogue of Theorem \ref{T:singleton}, where we denote unitary equivalence classes and the quotient operator topology by $[\cdot]_u$ and $QOT_u$, respectively.

\begin{corollary} \label{T:singleton2}
Let $p$ be a projection in a von Neumann algebra $\M$.  Then $\overline{\U(p)}^s$ consists of projections.
\begin{enumerate}
\item If $\M$ is finite,
$$\overline{\{[p]_u\}}^{QOT_u} = \{[p]_u\}.$$
\item If $\M$ is properly infinite and $p$ is finite, $$\overline{\{[p]_u\}}^{QOT_u} = \{[q]_u \mid [q]_u \leq [p]_u\}.$$
\item If $\M$ is properly infinite and $p^\perp$ is finite, $$\overline{\{[p]_u\}}^{QOT_u} = \{[q]_u \mid [q]_u \geq [p]_u\}.$$
\item If $p$ and $p^\perp$ are properly infinite and $c(p) = c(p^\perp)=1$, $$\overline{\{[p]_u\}}^{QOT_u} = (\p(\M)/\sim_u).$$
\end{enumerate}
\end{corollary}

\begin{proof} That $\p(\M)$ is strongly closed is a standard fact which may be proved in many ways, including Corollary \ref{T:spinc}.  The equalities will follow readily from Theorems \ref{T:main}(1) and \ref{T:singleton}; we simply need to describe the projections $q$ such that
$$q \in p, \qquad q^\perp \in p^\perp.$$

For $\M$ finite, ``$\in$" is Murray-von Neumann equivalence, so $q \sim_u p$.  In the rest of this paragraph assume $\M$ properly infinite.  When $p$ and $(1-p)$ are properly infinite with $c(p) = c(1-p) = 1$, there are no restrictions on $q$.  When $p$ is finite, we have the single restriction $q \preccurlyeq p$.  Say that $q \sim p_0 \leq p$.  Since the complement of a finite projection is Murray-von Neumann equivalent to the identity in a properly infinite algebra, $(1-q)$ and $(1-p_0)$ are Murray-von Neumann equivalent to each other.  It follows that $q \sim_u p_0 \leq p$ as desired.  The case where $(1-p)$ is finite is similar.
\end{proof}

At first glance this might have been unexpected: Murray-von Neumann equivalence, and not unitary equivalence of projections, is the right tool in an investigation of closed unitary orbits.  The statement $q \in p$ is strictly weaker than the analogous $q \in_u p$ (meaning $[q]_u \in \overline{\{[p]_u\}}^{QOT_u}$), and ``$\in_u$" cannot be substituted for ``$\in$" on the right-hand side of \eqref{E:main}.  The example in Remark \ref{T:noclosed} already illustrates this, as $\chi_{(0,2)}(k) \notin_u 1 = \chi_{(0,2)}(h)$.

\section{Crude and cruder: some reformulations} \label{S:reform}

\begin{definition}
Let $h$ be a normal operator in the von Neumann algebra $\M$.  The \textbf{crude $\M$-multiplicity function} for $h$ is the map
$$M_h: \{\text{open subsets of }\C \} \to \ps, \qquad \Os \mapsto [\chi_\Os(h)].$$
\end{definition}
Since the ambient algebra will always be clear in this paper, we may safely eschew a more explicit and cumbersome notation like ``$M^\M_h$."  This comment applies to Definition \ref{T:cruderdef}(2) as well.

\smallskip

From Theorem \ref{T:main}(2), for normal $h$ and $k$ we have $\overline{\U(h)}^{\|\|} = \overline{\U(k)}^{\|\|}$ if and only if $M_h = M_k$.  But in general we cannot use crude $\M$-multiplicity functions to describe strong* closures of unitary orbits.  Informally, this is because Murray-von Neumann equivalence distinguishes between different sizes of infinity, and the relation ``$\in$" does not.  So we introduce an equivalence relation which is similarly apathetic.

\begin{definition} \label{T:cruderdef} $\quad$
\begin{enumerate}
\item We define the \textbf{cruder} equivalence relation $\sim_c$ on $\p(\M)$ by $p \sim_c q$ iff $T(p) = T(q)$ for some (hence any) extended center-valued trace $T$.  We denote equivalence classes by $[\cdot]_c$.
\item For a normal operator $h \in \M$, we define the \textbf{cruder $\M$-multiplicity function} to be the map
$$M_h^c: \{\text{open subsets of }\C \} \to (\p(\M)/\sim_c), \qquad \Os \mapsto [\chi_\Os(h)]_c.$$
\end{enumerate}
\end{definition}

From \eqref{E:finclos} and \eqref{E:infclos3}, $p \sim_c q$ is the same as requiring both $p \in q$ and $q \in p$, and this in turn is the same as $\overline{\{[p]\}} = \overline{\{[q]\}}$ in ($\ps$, $QOT$).  This shows that $(\p(\M)/\sim_c)$ can be constructed as the maximal topological quotient of $(\ps, QOT)$ for which the quotient topology is $T_0$.  By \cite[Proposition 3.7]{S2005} or \eqref{E:par}, cruder equivalence agrees with Murray-von Neumann equivalence on $\p(\M)$ if and only if $\M$ is a direct sum of $\sigma$-finite von Neumann algebras.  And $(\p(\M)/\sim_c)$ clearly carries a (quotient) partial order given by
$$[p]_c \leq [q]_c \iff T(p) \leq T(q), \qquad p,q \in \p(\M).$$
(In fact $(\p(\M)/\sim_c)$ can be shown to be a complete lattice, via \cite[Theorem 6.1]{S2005} or more direct methods.)

Still another way to describe cruder equivalence is the following.  Let $z\M$ be the finite summand of $\M$.  Let $r$ be a projection with full central support such that $zr = z$ and $(1-z)r$ is a sum of centrally orthogonal $\sigma$-finite projections.  Then
\begin{equation} \label{E:cunit}
[p]_c \leq [q]_c \quad \iff \quad ([p] \wedge [r]) \leq ([q] \wedge [r]), \qquad p,q \in \p(\M).
\end{equation}
With \eqref{E:cunit} in mind, the reader will recognize the next result as the general von Neumann algebra version of Theorem \ref{T:strongbh}, with the class $[r]$ playing the role of $\aleph_0$.

\begin{theorem} \label{T:crudedom}
Let $h,k \in \M$ be normal.  Then
\begin{equation} \label{E:crudedom}
k \in \overline{\U(h)}^s \quad \iff \quad M^c_k \leq M^c_h.
\end{equation}
\end{theorem}

\begin{proof}
By considering central summands, it is enough to treat separately the cases of finite and properly infinite $\M$.  For $\M$ properly infinite, Theorem \ref{T:main}(1) and \eqref{E:infclos3} imply the conclusion.

When $\M$ is finite, Theorem \ref{T:main}(1) shows that $k \in \overline{\U(h)}^s$ if and only if $M^c_k = M^c_h$.  It suffices to prove that $M^c_k \leq M^c_h$ implies $M^c_k = M^c_h$.  Let $\Os \subset \C$ be an arbitrary open set.  For each $n$, consider the open set $\U_n = \Os^c + B_{1/n}(0)$.  We have that $\Os^c = \cap \U_n$, and so
\begin{align*}
T(\chi_{\Os^c}(k)) &= T(w-\lim \chi_{U_n}(k)) = w-\lim T(\chi_{U_n}(k)) \\ &\leq w-\lim T(\chi_{U_n}(h)) = T(w-\lim \chi_{U_n}(h)) = T(\chi_{\Os^c}(h)).
\end{align*}
Subtracting this from the identity operator gives $T(\chi_\Os(k)) \geq T(\chi_\Os(h))$.  But the opposite inequality holds by hypothesis, so we must have equality as desired. 
\end{proof}

Below we summarize our main progress in these terms.  All of these statements are direct consequences of previous results.

\begin{theorem} \label{T:reform} Let $h$ and $k$ be normal elements in the von Neumann algebra $\M$.
\begin{enumerate}
\item We have $k \in \overline{\U(h)}^{\|\|}$ if and only if $M_k = M_h$; the crude $\M$-multiplicity function is a complete invariant for norm-closed unitary orbits of normal operators in $\M$.
\item We have $k \in \overline{\U(h)}^s$ if and only if $M^c_k \leq M^c_h$; the cruder $\M$-multiplicity function is a complete invariant for strong-closed (also strong*-closed) unitary orbits of normal operators in $\M$.
\item Cruder $\M$-multiplicity functions can be identified with crude $\M$-multiplicity functions exactly when $\M$ is a direct sum of $\sigma$-finite algebras.
\item (combining (1), (2), and (3) above) If $\M$ is a direct sum of $\sigma$-finite algebras,
$$\overline{\U(h)}^{\|\|} = \overline{\U(k)}^{\|\|} \iff \overline{\U(h)}^{s^*} = \overline{\U(k)}^{s^*} \iff \overline{\U(h)}^s = \overline{\U(k)}^s.$$
\end{enumerate}
\end{theorem}

Ding and Hadwin have noted the agreement of approximate equivalence and weak approximate equivalence for *-homomorphisms from certain $C^*$-algebras to von Neumann algebras with separable predual (\cite[Corollary 11]{DH}).  This statement, which we do not make precise here, contains Theorem \ref{T:reform}(4) for $\M$ with separable predual.  (We mentioned in the introduction that their work also subsumes Theorem \ref{T:reform}(1) when $\M$ has separable predual.)

For $\M = \B(\h)$, all of Theorem \ref{T:reform} follows from Hadwin's results about ``approximate multiplicity" (\cite[Propositions 5.3, 5.4 and Theorem 7.1]{H1981}).  This is a noncommutative version of crude multiplicity which generalizes points in the spectrum to irreducible representations or primitive ideals.  (Hadwin's analogue for the cruder multiplicity function appears in his Propositions 5.3(2) and 5.4(4).)  Underlying the usefulness of approximate multiplicity -- and it is a key ingredient in the proof of Theorem \ref{T:hv} -- is a special feature of the dimension theory of $\B(\h)$ which we now discuss.

\begin{remark} \label{T:ptcrude}
There is an important difference between $\B(\h)$ and general $\M$ which has influenced our approach in this paper: $\ps$ is not necessarily well-ordered.  Indeed, it follows from dimension theory (\cite{To} or \cite[Corollary 2.8]{S2005}) that $\ps$ is well-ordered if and only if $\M$ is a factor of type I or III.

It is a daily chore of operator theory to interpret the multiplicity of points in the spectrum of an operator.  Our sources \cite{EEL,AD,D1986} all employ (somewhat interdependently) the cardinal
\begin{equation} \label{E:ptcrude}
\inf_{r>0} \text{rank} (\chi_{B_r(\lambda)}(h)),
\end{equation}
where $h$ is a normal operator in some $\B(\h)$.  This generates a ``crude multiplicity function" whose domain is $\C$.  It contains all the same information as the crude multiplicity function that we have defined, and it simplifies certain kinds of arguments by offering ready-made discrete approximants to the spectral measure (\cite{AD,D1986}).  One can mimic \eqref{E:ptcrude} in general von Neumann algebras by replacing rank with equivalence class, but this tends to lose too much information for our purposes, exactly because $(\p(\M)/\sim)$ is not well-ordered.  For normal $h$ with no eigenprojections in a $\text{II}_1$ factor, the quantity would be identically zero.  Nonetheless several authors have used a von Neumann algebra version of \eqref{E:ptcrude} to pick out the essential spectrum: \cite[Corollary 2.3]{F}, \cite[Proposition 3.8]{Kaf1} (see Corollary \ref{T:allin} below), \cite[Lemma 1.2]{Hi} (even sorting the essential spectrum by weight). 

\end{remark}

\section{Essential spectrum} \label{S:esssp}

Let $\M$ be a properly infinite von Neumann algebra, and let $\kom$ be the norm-closed two-sided ideal generated by the finite projections.  (The symbol $\kom$ only has this meaning in the present paper.)  An element of $\M$ is called \textit{essentially normal} if its image in the quotient $C^*$-algebra $\M/\kom$ is normal.  The \textit{(left, right) essential spectrum} of $h \in \M$, denoted $\text{sp}^e(h)$ (resp. $\text{sp}^{le}(h)$, $\text{sp}^{re}(h)$), is the (resp. left, right) spectrum of $h + \kom$ in $\M/\kom$.  (The left/right spectrum of an element $x$ in a unital Banach algebra is the set of complex scalars $\lambda$ such that $x - \lambda I$ fails to be left/right invertible.)

Next we review some ideas of Halpern.  We only need to define his terms for the ideal $\kom$, but they are meaningful for other central ideals (\cite{Halp1972}) as well.  It is appropriate to remark that Str\u{a}til\u{a} and Zsid\'{o} (\cite{SZ1973, St1973}) developed similar machinery at the same time as Halpern.  

\begin{definition}(\cite[Definition 3.6]{Halp1972}) \label{D:ecs}
Let $h$ be an operator in a properly infinite von Neumann algebra $\M$.  Recall that elements of $\z(\M)$ are identified by the Gelfand transform with continuous functions on the maximal ideal space $\ozm$.  For $\zeta \in \ozm$, $\zeta + \kom$ means the sum ideal in $\M$.

The \textbf{(left, right) essential central spectrum} of $h$, written $\text{sp}^{ec}(h)$ (resp. $\text{sp}^{lec}(h)$, $\text{sp}^{rec}(h)$), is the set of central elements $z$ such that $h - z + (\zeta + \kom)$ is not (resp. left, right) invertible in $\M/(\zeta + \kom)$ for any $\zeta \in \ozm$.
\end{definition}

We inject a few obvious comments.  If $\M$ is a factor, the essential central spectra can be identified with the (numerical) essential spectra, as multiples of the identity.  And if $h$ is essentially normal, all three essential central spectra agree.

Less obviously, the three essential central spectra are never empty.

\begin{theorem}$($\cite[Theorem 3.5]{Halp1972}$)$ \label{T:nonempty}
For any operator $h$ in a properly infinite von Neumann algebra $\M$, the set $\mbox{\textnormal{sp}}^{lec}(h) \cap \mbox{\textnormal{sp}}^{rec}(h)$ is nonempty.
\end{theorem}

\begin{proposition} $(\sim$\cite{Halp1972}$)$ \label{T:halp}
For a self-adjoint operator $h$ in a properly infinite von Neumann algebra $\M$, the following conditions are equivalent:
\begin{enumerate}
\item $0 \in \mbox{\textnormal{sp}}^{ec}(h)$;
\item there is a sequence $\{p_n\}$ of pairwise orthogonal properly infinite projections with full central support satisfying $\|h p_n\| < 1/n$;
\item for any $\e > 0$, $\chi_{(-\e, \e)}(h)$ is properly infinite with full central support.
\end{enumerate}
\end{proposition}

\begin{proof} The implication (1) $\to$ (3) follows from \cite[Proposition 3.13]{Halp1972}; (1) $\leftrightarrow$ (2) is \cite[Corollary 3.16]{Halp1972}.  To see (3) $\to$ (2), construct $p_n$ inductively by halving $(\chi_{(-1/n,1/n)}(h) \wedge (\sum_{j=1}^{n-1} p_j)^\perp)$ (omitting the sum for $n=1$).
\end{proof}

For use in Section \ref{S:agree}, we record the following corollary.  A non-factor version was obtained by Kaftal (\cite[Proposition 3.8]{Kaf1}) without direct appeal to Halpern's theory.

\begin{corollary} \label{T:allin}
For $\lambda \in \C$ and normal $h$ in an infinite factor $\M$,
$$\lambda \in \mbox{\textnormal{sp}}^e(h) \iff \chi_\Os(h) \mbox{\textnormal{ is infinite, }} \forall \mbox{\textnormal{ open }} \Os \ni \lambda.$$
\end{corollary}

\begin{proof}
Using Proposition \ref{T:halp} at the key step,
\begin{align*} \lambda \in \text{sp}^e(h) &\iff \lambda I \in \text{sp}^{ec}(h) \\ &\iff 0 \in \text{sp}^{ec}(h - \lambda I) \\ &\iff 0 \in \text{sp}^{ec}(|h - \lambda I|) \\ &\iff \forall \e > 0, \; \chi_{(-\e, \e)}(|h - \lambda I|) \text{ is infinite} \\ &\iff \forall \e > 0, \; \chi_{B_\e(\lambda)}(h) \text{ is infinite}. \qedhere
\end{align*}
\end{proof}

Here is the main result of the section.

\begin{theorem} \label{T:lec}
If $h$ is any operator (not necessarily normal) in a properly infinite von Neumann algebra $\M$, then
\begin{equation} \label{E:lec}
\overline{\U(h)}^s \cap \z(\M)  = \mbox{\textup{sp}}^{lec}(h).
\end{equation}
\end{theorem}

\begin{proof} Let $z \in \z(\M)$.  We will argue by a chain of logical equivalences.
\begin{align*}
\exists \{u_\alpha\} &\subset \U(\M), \quad u_\alpha h u_\alpha^* \overset{s}{\to} z \\ &\iff \exists \{u_\alpha\} \subset \U(\M), \quad u_\alpha (h - z) u_\alpha^* \overset{s}{\to} 0 \\ &\iff \exists \{u_\alpha\} \subset \U(\M), \quad u_\alpha |h - z| u_\alpha^* \overset{s}{\to} 0.
\end{align*}
Theorem \ref{T:main}(1) tells us that the last condition is equivalent to 
$$\forall \e > 0, \quad \chi_{[0,\e)}(|h - z|) \text{ is properly infinite with full central support}.$$
By Proposition \ref{T:halp}, this is equivalent to $0 \in \text{sp}^{ec}(|h - z|).$  From elementary invertibility considerations the latter is equivalent to $0 \in \text{sp}^{lec}(h-z)$ and finally $z \in \text{sp}^{lec}(h)$.
\end{proof}

What about the strong* closure of a unitary orbit?  It is easy to see that for $h$ in a properly infinite von Neumann algebra $\M$,
\begin{align*} \overline{\U(h)}^{s^*} \cap \z(\M) &\subseteq \overline{\U(h)}^s \cap \left( \overline{\U(h^*)}^s \right)^* \cap \z(\M) \\ &= \text{sp}^{lec}(h) \cap \text{sp}^{rec}(h),
\end{align*}
and the last set is nonempty, by Corollary \ref{T:nonempty}.  But in general the inclusion may be proper, as we now demonstrate.

\begin{corollary} \label{T:dixmier}
In a properly infinite von Neumann algebra, the strong-closed unitary orbit of any operator meets the center.

This becomes false when ``strong" is replaced by ``strong*."
\end{corollary}

\begin{proof}
The first statement follows from Theorems \ref{T:nonempty} and \ref{T:lec}.

Now take $\M = \B(\ell^2)$, and let $v$ be a partial isometry between complementary projections.  If $u_\alpha v u_\alpha^* \overset{s^*}{\to} w$, then
$$w^* w = s^*-\lim u_\alpha (v^* v) u_\alpha^* \in \p(\M);$$
$$w^*w + w w^* = s^*-\lim u_\alpha (v^* v + v v^*) u_\alpha^* = 1.$$
So $w$ is also a partial isometry between complementary projections, and therefore $w$ is unitarily conjugate to $v$.  This shows $\overline{\U(v)}^{s^*} = \U(v)$, which clearly does not intersect the center at all.

This phenomenon is directly attributable to the fact that $C^*(v+\kom) \subseteq \B(\ell^2)/\kom$ has no one-dimensional representations (\cite{PS}).
\end{proof}

Hadwin's work on operator-valued spectra (\cite{H1977,H1987}, see also \cite{E} and \cite{BD}) is relevant here.  Specializing to $\B(\ell^2)$ and omitting definitions, Corollary \ref{T:dixmier} is the observation that the one-dimensional essential reducing operator spectrum may be empty, while the one-dimensional left essential operator spectrum never is.

\smallskip

Corollary \ref{T:dixmier} is related to the topic of ``Dixmier averaging," which descends from Dixmier's original theorem (\cite{Di1949}) that the norm-closed convex hull of a unitary orbit always intersects the center.  Broadly defined, this industry is concerned with the following question: if a subgroup of $\U(\M)$ acts by conjugation on a subset of $\M$, must there be a fixed point?  A usual technique is to average, using convexity of the subset, but there is no convexity (or compactness) in Corollary \ref{T:dixmier}.  (As we demonstrate in a forthcoming paper with C. Akemann, convexity sometimes appears automatically if one considers weak-closed unitary orbits instead.  A prototype result is \cite[Corollary 2.3]{Halp1977}.) 

Zsid\'{o} asked (\cite[Problem 2]{Z1974}) whether the norm-closed unitary orbit of an operator in a properly infinite algebra necessarily meets the center.  We suppose that a different topology was intended, as the answer is obviously no.  (A norm-closed unitary orbit which intersects the center is necessarily a singleton in the center, cf. Theorem \ref{T:strongnorm}.)  In any case Corollary \ref{T:dixmier} provides a ``yes" answer for topologies as coarse as the strong and a ``no" answer for topologies as fine as the strong*.  

\smallskip

For comparison with Theorem \ref{T:lec}, we give the following

\begin{proposition}
For an operator $h$ in a finite von Neumann algebra $\M$, $$\overline{\U(h)}^s \cap \z(\M) \neq \varnothing \iff h \in \z(\M).$$
\end{proposition}

\begin{proof}
We need only to prove the forward implication.  Let $T$ be the center-valued trace.  If $u_\alpha h u_\alpha^* \overset{s}{\to} z \in \z(\M)$, then
\begin{align*} &(u_\alpha h u_\alpha^* - z)^*(u_\alpha h u_\alpha^* - z) \overset{s}{\to} 0 \\ &\Rightarrow 0 = s-\lim T[(u_\alpha h u_\alpha^* - z)^*(u_\alpha h u_\alpha^* - z) ] = T[(h-z)^*(h-z)] \\ &\Rightarrow h=z.
\qedhere
\end{align*}
\end{proof}

\smallskip

There is a ``picture" which goes with Theorem \ref{T:lec}.  In an arbitrary von Neumann algebra $\M$, consider the equivalence relation on normal operators defined by
$$h \sim k \iff \overline{\U(h)}^s = \overline{\U(k)}^s.$$
On equivalence classes we have a partial order defined by containment of strong-closed unitary orbits:
$$[k] \leq [h] \iff \overline{\U(k)}^s \subseteq \overline{\U(h)}^s \iff k \in \overline{\U(h)}^{s} \iff M^c_k \leq M^c_h.$$
If $\M$ is finite, this partial order is just equality.

If $\M$ is properly infinite, it follows from Theorem \ref{T:lec} that the minimal elements in this partial order are the singleton classes of central operators in $\M$.  In a complete atomic Boolean algebra, an element is completely determined by the minimal elements under it; in the situation at hand one recovers the essential central spectrum of a representative of the equivalence class, which typically does not suffice to recover the class itself.  We note the exception in the following theorem, whose title will be explained in the next section.

\begin{theorem} (Type III Weyl-von Neumann-Berg theorem) \label{T:wvnb'}
Let $h,k \in \M$ be normal operators in a type III von Neumann algebra.  Then
\begin{equation} \label{E:ecs3}
k \in \overline{\U(h)}^s \iff \mbox{\textnormal{sp}}^{ec}(k) \subseteq \mbox{\textnormal{sp}}^{ec}(h).
\end{equation}
\end{theorem}

\begin{proof}
The forward implication is clear from Theorem \ref{T:lec}.

Now assume $k \notin \overline{\U(h)}^s$.  Since projections in type III algebras are cruder equivalent to their central supports, by Theorem \ref{T:crudedom} there is an open set $\Os$ with $c(\chi_\Os(k)) \nleq c(\chi_\Os(h))$.  Setting
$$\Os_n = \{\lambda \mid \text{dist}(\lambda, \Os^c) > 1/n\},$$
we have $\chi_{\Os_n}(k) \nearrow \chi_\Os(k)$.  So there must be some $n$ for which  
$$c(\chi_{\Os_n}(k)) \nleq c(\chi_\Os(h)).$$
Set
\begin{equation} \label{E:restr}
z = \left( c(\chi_{\Os_n}(k)) \wedge (c(\chi_\Os(h)))^\perp \right) \neq 0, \qquad p = z\chi_{\Os_n}(k).
\end{equation}
We view $z\M$ as the ambient algebra until the last paragraph of the proof.  Write $zk = pk + p^\perp k$ and apply Corollary \ref{T:dixmier} to $pk \in p\M p$:
$$\exists y_1 \in \overline{ \{ upku^* \mid u \in \U(p\M p) \} }^s \cap \z(p\M p) \quad \Rightarrow \quad (y_1 + p^\perp k) \in \overline{\U(zk)}^s.$$

Since $c(p) = z$, multiplication by $p$ gives an isomorphism from $\z(z\M)$ to $\z(p \M p)$ (\cite[Lemma 6.39]{AS}).  Use this to write $y_1 = p y$, with $y \in \z(z\M)$.  We make two claims about $y$:
\begin{itemize}
\item $y \in \overline{\U(zk)}^s$;
\item $y \notin \overline{\U(zh)}^s$.
\end{itemize}

Consider the first claim, and keep in mind that $y$ is central.  For any open $U$, $\chi_U(y_1) = p \chi_U(y)$ in $p\M p$.  Taking central supports in the larger algebra $z\M$,
$$c(\chi_U(y_1)) = c(p \chi_U(y)) = c(p) \chi_U(y) = \chi_U(y).$$
Therefore
$$[\chi_U(y)]_c = [\chi_U(y_1)]_c \leq [\chi_U(y_1) + \chi_U(p^\perp k)]_c = [\chi_U(y_1 + p^\perp k)]_c \leq [\chi_U(zk)]_c.$$
We used Theorem \ref{T:crudedom} for the second inequality, and we may use it again to conclude that $y \in \overline{\U(zk)}^s$.

For the second claim, the isomorphism and \eqref{E:spinc2} imply
$$\text{sp}_{z\M}(y) = \text{sp}_{p\M p}(y_1) \subseteq \text{sp}_{p\M p}(pk) \subseteq \overline{\Os_n} \subseteq \Os.$$
(Here the subscripts simply clarify the ambient algebra.)  From \eqref{E:restr}, 
$$\text{sp}_{z\M}(zh) \cap \Os = \varnothing,$$
so by \eqref{E:spinc2} again $y \notin \overline{\U(zh)}^s$.

The two claims show $\text{sp}^{ec}(zk) \nsubseteq \text{sp}^{ec}(zh)$ in $z\M$.  Since taking the essential central spectrum commutes with restriction to a central summand (\cite[Proposition 3.10]{Halp1972}), we have $\text{sp}^{ec}(k) \nsubseteq \text{sp}^{ec}(h)$, finishing the proof.
\end{proof}

\begin{remark} \label{T:sz}
In a type III algebra, $\kom = \{0\}$, so the essential central spectrum of any operator $x$ is just the ``central spectrum."  We have not defined this term separately in the present article, but it is explored in \cite[Section 6]{SZ1973} and usefully characterized as
$$\{z \in \z(\M) \mid (x - z)y \text{ is not invertible in $y\M$ for any }y \in \p(\z(\M))\}.$$
By a variant of (\cite[Lemma 1.3]{St1973}), this kind of description also holds for (left, right) essential spectra with respect to $\kom$ or other central ideals. 
\end{remark}

\section{Some relations to the literature} \label{S:rel}
\subsection{The Weyl-von Neumann-Berg theorem in von Neumann algebras} \label{SS:wvnb} This theorem is really a body of results due to several authors; we choose this title (WvNB) for its currency and relative brevity.  For us, the main statement is this: \textit{Two normal operators in $\B(\ell^2)$ are unitarily equivalent modulo the compact operators if and only if they have the same essential spectrum.}  Primary sources for this version are \cite{W1909,vN,Be,S1970,Hal1972,S1971,EEL}.   A fuller history than we can offer here would discuss unbounded operators, nonseparable Hilbert spaces, Hilbert-Schmidt differences, other extensions of Weyl's original theorem, and most notably, the beautiful generalizations to representations found by Voiculescu (\cite{V1976}) and Hadwin (\cite{H1981}, cf. Theorem \ref{T:hv} of this paper).  We instead draw attention to some von Neumann algebra versions of WvNB.

Let $h,k \in \M$ be normal.  When $\M$ is a $\sigma$-finite semifinite infinite factor, a WvNB theorem (\cite{Z1975,F,Kaf2,K1988}) gives the equivalence of
\begin{enumerate}
\item $h$ and $k$ have the same essential spectrum;
\item $h$ and $k$ are unitarily equivalent mod $\kom$ -- i.e., there exist $c \in \kom$ and $u \in \U(\M)$ such that $k = uhu^* + c$.
\end{enumerate}
A version of this theorem has also been proved for properly infinite semifinite non-factors with separable predual (\cite{CK}), via disintegration over the center.

If $\M$ is $\B(\ell^2)$, one also has the equivalence of the stronger conditions
\begin{enumerate}
\item[(1')] $h$ and $k$ have the same strong*-closed unitary orbit;
\item[(2')] $h$ and $k$ have the same norm-closed unitary orbit;
\item[(3')] same as (2), except $c$ can be chosen with arbitrarily small norm.
\end{enumerate}
In appropriate language, this was actually proved for representations of (noncommutative) $C^*$-algebras (\cite{V1976,A1977,H1981}).  The jump between the two sets of conditions was spurred by Brown-Douglas-Fillmore theory (\cite{BDF}), which considered condition (2) for essentially normal operators.

What about type III algebras?  If $\M$ is a factor, \eqref{E:spinc2} and Corollary \ref{T:converse}(2) assert the equivalence of conditions (1) and (1').  The equality case in Theorem \ref{T:wvnb'} is the correct extension of this WvNB-type result to non-factors.  If moreover the algebra is a direct sum of $\sigma$-finite algebras, (1') and (2') are also equivalent by Theorem \ref{T:reform}(4).  (Conditions (2) and (3') reduce to unitary equivalence in a type III algebra, since $\kom = \{0\}$.)

\subsection{Distance between unitary orbits} \label{SS:dist}
For normal $h$ and $k$ in a von Neumann algebra $\M$, define the \textit{spectral distance} between them as
\begin{align} \label{E:spdist1} &\delta (h,k) = \mbox{\textnormal{infimum of the nonnegative reals $r$ satisfying}}\\ \notag [&\chi_{\Os + B_r(0)}(h)] \geq [\chi_\Os(k)], \quad [\chi_{\Os + B_r(0)}(k)] \geq [\chi_\Os(h)], \quad \forall \mbox{\textnormal{ open }} \Os.
\end{align}
This kind of definition was first proposed for $\B(\h)$ in \cite{AD}, and a trace version is stated for $\sigma$-finite semifinite factors in \cite{HN}.  If $\M$ is the matrix algebra $M_n$, and $h$ and $k$ have eigenvalue lists $\{\lambda_j\}$ and $\{\mu_j\}$, respectively, then $\delta(h,k)$ works out to be
\begin{equation} \label{E:spdist2}
\min_\sigma \max_j |\lambda_j - \mu_{\sigma(j)}|,
\end{equation}
where $\sigma$ runs over the permutations on $n$ letters.

Spectral distance has been much studied in connection with the famous problem of determining the distance between the unitary orbits of two normal operators.  Weyl (\cite{W1912}) showed in 1912 that \eqref{E:spdist2} gives the answer when $h$ and $k$ are self-adjoint matrices; since at least the 1960s this was conjectured for normal matrices as well.  There were several developments in the 1980s, including the solution \eqref{E:spdist1} for (infinite-dimensional) self-adjoint operators (\cite{AD}), an inequality for normal operators (\cite{BDM,D1986,BDK}), and a von Neumann algebra version of this inequality (\cite{HN}), which we now state. 
 
\begin{theorem}$($\cite{HN}$)$ \label{T:hn}
Let $h,k$ be normal elements in a $\sigma$-finite factor $\M$.  Then
\begin{equation} \label{E:hn}
\mbox{\textnormal{dist}}(\U(h), \U(k)) \leq \delta(h,k) \leq c \:\mbox{\textnormal{dist}}(\U(h), \U(k)),
\end{equation}
where $c$ is a universal constant known to be $< 2.91$.
\end{theorem}
The conjecture was that $c=1$, but in 1992 it was shown by Holbrook (\cite{Ho}) that $c>1$ already for $3 \times 3$ matrices.  For more background on this problem, see the survey \cite{DS}.

An immediate consequence of \eqref{E:hn} is
\begin{equation} \label{E:sdelta}
k \in \overline{\U(h)}^{\|\|} \iff \delta(h,k) = 0.
\end{equation}
We make an observation which shows how \eqref{E:sdelta} implies Theorem \ref{T:main}(2) (for $\sigma$-finite factors).

\begin{proposition} \label{T:delta}
Let $h$ and $k$ be normal elements in a $\sigma$-finite factor $\M$.  Then
\begin{equation} \label{E:delta}
\delta(h,k) = 0 \iff [\chi_\Os(h)] = [\chi_\Os(k)], \quad \forall \text{ open } \Os \subset \C.
\end{equation}
\end{proposition} 

\begin{proof}
Suppose $\delta(h,k)=0$, and let $\Os \subset \C$ be open.  Define
$$\Os_n = \{\lambda \in \C \mid \text{dist}(\lambda, \Os^c) > 1/n\}, \qquad n = 1,2,\dots$$
Note that
$$\Os_n + B_{1/n}(0) \subset \Os \quad \Rightarrow \quad [\chi_{\Os_n}(k)] \leq [\chi_{\Os_n + B_{1/n}(0)}(h)] \leq [\chi_\Os(h)].$$

By Proposition \ref{T:normal},
$$[\chi_\Os(k)] = [\chi_{\cup \Os_n}(k)] = [\vee \chi_{\Os_n}(k)] = \vee [\chi_{\Os_n}(k)] \leq [\chi_\Os(h)].$$
Since the roles of $h$ and $k$ may be interchanged, we have $[\chi_\Os(h)] = [\chi_\Os(k)]$ as desired.

The reverse implication is trivial.  
\end{proof}

\begin{remark} Proposition \ref{T:delta} remains true when $\M$ is an arbitrary von Neumann algebra; the fact is that $[\vee p_n] = \vee [p_n]$ always holds for increasing \textit{sequences} of projections.  Here is a one-sentence outline of the case not covered by Proposition \ref{T:normal}, assuming that $\vee [p_n]$ is a properly infinite class.  Let $q_1 = p_1$ and $q_n = p_n - p_{n-1}$ for $n \geq 2$, rewrite the join as a sum, and apply complete additivity: 
$$[ \vee p_n ] = \left[ \sum q_n \right] = \sum [q_n] \leq \aleph_0 \cdot (\vee [p_n]) = \vee [p_n] \leq [ \vee p_n].$$

\end{remark}  

Clearly \eqref{E:sdelta} and \eqref{E:delta} entail \eqref{E:main2} (for $\sigma$-finite factors). 

In the present paper we only characterize when two norm-closed unitary orbits are equal.  There is no direct contribution here to the solution of the distance problem and its variants, but it seems possible that Theorem \ref{T:main}(2) and Proposition \ref{T:delta} could be useful in extending a solution, or an estimate like \eqref{E:hn}, to all von Neumann algebras.

\section{Agreements between $\U(h)$ and its closures} \label{S:agree}

For any operator $h$ in a von Neumann algebra $\M$, we obviously have the inclusions
$$\U(h) \subseteq \overline{\U(h)}^{\|\|} \subseteq \overline{\U(h)}^{s^*} \subseteq \overline{\U(h)}^s.$$
The main goal of this section is to characterize when each of these inclusions becomes an equality, under the assumption that $h$ is normal.  This is done in Theorems \ref{T:strongnorm}, \ref{T:ss}, and \ref{T:normc}.  As must be expected, there is quite a bit of relevant literature for the case $\M = \B(\ell^2)$.

We begin this section by explaining some of the difficulties involved in the description of the unclosed orbit $\U(h)$.  These are hardly new ideas, but they serve to contrast with our main results and facilitate the proof of Theorem \ref{T:normc}.

\smallskip

We have already mentioned that the problem of unitary equivalence was solved for $\B(\h)$ by the hands of several authors.  Probably the best-known work is the 1951 book of Halmos (\cite{Hal1951}), which is devoted to a succinct answer: ``Two normal operators $h$ and $k$ in $\B(\h)$ are unitarily equivalent if and only if their spectral measures have the same \textit{multiplicity function}."

In a willful oversimplification, one may isolate two main steps in this solution (assuming for simplicity that $\h$ is separable).  First, decompose $W^*(h)$ and $W^*(k)$ as direct sums of multiples of maximal abelian *-subalgebras (heretofore called ``MASAs") on subspaces of $\h$.  By comparing the operator summands of equal multiplicity, this reduces the problem to one in which both operators generate a MASA.  Second, solve the latter problem by a spatial spectral theorem: such an operator is unitarily equivalent to multiplication by the function $f(z) = z$ on some $L^2(\C, \mu)$, where the possibly infinite measure $\mu$ is only canonical up to mutual absolute continuity.  And this is enough, as two such multiplication operators are unitarily equivalent exactly when the measures are mutually absolutely continuous.

This version of multiplicity theory more or less amounts to identifying unitary equivalence classes of representations of singly-generated abelian $C^*$-algebras.  It is a small jump to remove ``singly-generated;" more interesting (and less successful) are the attempts to remove ``abel\-ian."  In the previous paragraph, the decomposition according to multiplicity is exactly the central decomposition of the type I commutant into $\text{I}_n$ summands.  This suggests, accurately, that a multiplicity theory for algebra representations will go well as long as all commutants are type I.  Such algebras and operators have many names in the literature, and their multiplicity theory is accessibly presented in Arveson (\cite{A1976}).  For a brave foray away from type I to arbitrary operators, see Ernest (\cite{E}).  A more abstract approach to multiplicity is given by Kadison in \cite{K1957}. 

The two-step program above has other difficulties when applied to the problem of unitary equivalence for normal operators in a general von Neumann algebra $\M$.  Regarding multiplicity, the \textit{relative} commutant of a normal operator need not be type I.  See Bures (\cite{Bu}) for work along these lines (and note that ``multiplicity" has, well, multiple meanings in the literature).  More fundamentally, two normal operators which generate MASAs and have identical spectral data need not be unitarily conjugate.  The following example was independently chosen to demonstrate a similar point in \cite{AK}.  The underlying mathematical ideas are mostly due to Dixmier (\cite{Di1954}).

\begin{example} \label{T:masas}
A MASA in a von Neumann factor $\A \subset \M$ is said to be \textit{regular} if its normalizer (=$\{u \in \U(\M) \mid u\A u^* = \A\}$) generates the algebra $\M$.  The hyperfinite $\text{II}_1$ factor $\RR$ contains both regular and non-regular MASAs.

Let $\A_1 \subset \RR$ be a regular MASA, $\A_2 \subset \RR$ a non-regular MASA.  We equip $\RR$ with its normalized trace $\tau$.  Each of $(\A_1, \tau)$ and $(\A_2, \tau)$ is *-isomorphic, in a state-preserving way, with $(L^\infty[0,1], m)$.  (See the proof of \cite[Theorem III.1.22]{T}; $m$ denotes integration against Lebesgue measure.)  Choosing such isomorphisms, let $h_j \in \A_j$ ($j=1,2$) be the positive operator which corresponds to the function $f(x) = x$ in $L^\infty[0,1]$.  We have that $\A_j = W^*(h_j)$.

The spectral data for $h_1$ and $h_2$ are indistinguishable: for any measurable set $E$ the corresponding spectral projections of $h_1$ and $h_2$ are Murray-von Neumann equivalent.  (Their traces both equal the Lebesgue measure of $E \cap [0,1]$.)  The ``multiplicities" of $W^*(h_j)$ are both one, since the algebras are maximal.  

But there is no automorphism of $\RR$ which takes $h_1$ to $h_2$, for any such would take $\A_1$ to $\A_2$, and the property of regularity is clearly invariant under automorphisms.

\end{example}

Example \ref{T:masas} shows that spectral data and multiplicity are too weak to distinguish automorphism orbits in a von Neumann algebra, much less unitary orbits.  There are just too many places for a MASA to sit.  On the other hand, the approximate equivalence of \eqref{E:main2} only requires some of the spectral data, and completely ignores the structure of the inclusion.  Herrero's book \cite{He} actually takes a general version of this phenomenon as a unifying theme for several types of operator approximation.

\smallskip

\textsc{The equality $\overline{\U(h)}^{\|\|} = \overline{\U(h)}^{s^*}$.} This is a fairly direct application of our work.

\begin{theorem} \label{T:strongnorm} Let $h$ be normal in a von Neumann algebra $\M$.
\begin{enumerate}
\item If $\M$ is finite, $\overline{\U(h)}^{\|\|}$ is strong-closed.
\item If $\M$ is properly infinite, the following conditions are equivalent:
\begin{enumerate}
\item[(a)] $\overline{\U(h)}^{\|\|}$ meets the center;
\item[(b)] $h$ belongs to the center;
\item[(c)] $\overline{\U(h)}^{\|\|} = \overline{\U(h)}^s$;
\item[(d)] $\overline{\U(h)}^{\|\|} = \overline{\U(h)}^{s^*}$.
\end{enumerate}
\end{enumerate}
\end{theorem}

\begin{proof}
The first statement follows from Theorem \ref{T:main}, Theorem \ref{T:singleton}(1), and the observation that the strong and strong* topologies agree in a finite algebra.  For self-adjoint elements in factors, this has been noted in \cite[Theorem 5.4]{AK}.

In the second statement, the implications (a) $\to$ (b) $\to$ (c) $\to$ (d) are obvious.  Assuming (d) and using the normality of $h$ and Corollary \ref{T:dixmier},
$$\overline{\U(h)}^{\|\|} \cap \z(\M) = \overline{\U(h)}^{s^*} \cap \z(\M) = \overline{\U(h)}^s \cap \z(\M) \neq \varnothing. \qedhere$$
\end{proof}

For non-normal elements, conditions (a), (b), (c) in Theorem \ref{T:strongnorm}(2) are equivalent by the same reasoning.  But (d) is strictly weaker, as is seen by consideration of the operator $v$ from the proof of the second statement in Corollary \ref{T:dixmier}.

\smallskip

\textsc{The equality $\overline{\U(h)}^{s^*} = \overline{\U(h)}^s$.} The proof of our characterization requires an array of tools from operator theory, so to mesh cleanly with existing results, we stipulate that the von Neumann algebras be factors with separable predual.

\begin{definition} Consider an operator $x$ on $\ell^2$.  We say $x$ is \textbf{reductive} if for any projection $p$, $xp = pxp$ implies $xp = px$.  We say $x$ is \textbf{strongly reductive} if for every $\e>0$ there is $\delta > 0$ such that any projection $p$ satisfying $\|(1-p)xp\| < \delta$ must also satisfy $\|xp - px\| < \e$.

In another vein, we say $x$ is \textbf{subnormal} if it has a normal extension, i.e. $x$ is the restriction of a normal operator to an invariant subspace.

Finally, we say $x$ is \textbf{pure} if it has no normal restriction to a reducing subspace.
\end{definition}

Here are most of the facts we need about subnormal operators.

\begin{enumerate}
\item The subnormal operators are precisely the strong limits of normal operators.
\item A subnormal operator $x$ is \textit{hyponormal}; i.e. $x^*x \geq x x^*$.
\item A subnormal operator $x$ has a \textit{minimal normal extension} which is unique up to unitary equivalence and denoted here by $\text{mne}(x)$. 
\item When $x$ is subnormal, $\text{sp}(\text{mne}(x))$ is obtained from $\text{sp}(x)$ by adjoining some of the bounded components of $\text{sp}(x)^c$.
\end{enumerate}
These can all be found in Chapter II of \cite{Co1991}, the standard reference for subnormal operators.

Next we prepare an easy lemma about subsets of the plane.  To avoid both needless repetition and pretensions of coining a new term, we will say that a compact subset of the complex plane is \textit{small} if it has connected complement and no interior.

\begin{lemma} \label{T:small} ${}$
\begin{enumerate} 
\item A subset of a small set is small.
\item A compact set with small boundary equals its boundary.
\end{enumerate}
\end{lemma}

\begin{proof}
For the first statement, let $A \subseteq B$, with $B$ small.  Then $A$ clearly has no interior.  If $A^c$ is disconnected, write $A^c = U_1 \cup U_2$, with $U_1, U_2$ disjoint, nonempty, and open.  Then
$$B^c \subseteq A^c \Rightarrow B^c = (B^c \cap U_1) \cup (B^c \cap U_2),$$
and by the smallness of $B$, one of these last two is empty.  But
$$B^c \cap U_j = \varnothing \Rightarrow U_j \subseteq B,$$
which contradicts smallness of $B$.

For any set $C$, we have the disjoint union
\begin{equation} \label{E:int}
(\partial C)^c = \text{int}(C) \cup \text{int}(C^c),
\end{equation}
where we write ``int" for ``interior of."  If $C$ is compact, then $\text{int}(C^c) = C^c$ and is nonempty.  If also $\partial C$ is small, by \eqref{E:int} $\text{int}(C)$ must be empty.  Now \eqref{E:int} reads $(\partial C)^c = C^c$, so $\partial C = C$.
\end{proof}

Soon after their introduction in the mid-1970s, strongly reductive operators were completely understood.

\begin{theorem} \label{T:har} $($\cite{Har, AFV}$)$ An operator is strongly reductive if and only if it is normal with small spectrum.
\end{theorem}

For operators on $\ell^2$, Hadwin gave useful descriptions of the strong and strong* closures of the unitary orbit as the sets of \textit{approximate restrictions} and \textit{approximate summands}, respectively (\cite[Theorem 4.4]{H1987}).  This means that for $x \in \B(\ell^2)$,
\begin{align}
\label{E:arest} &\overline{\U(x)}^s = \\ \notag &\{y \mid \exists \{v_n\}, v_n^* v_n = 1 \text{ and } \|v_n^*xv_n - y\|, \|(1- v_n v_n^*) x v_n v_n^* \| \to 0\}; \\
\label{E:asumm} &\overline{\U(x)}^{s^*} = \\ \notag &\{y \mid \exists \{v_n\}, v_n^* v_n = 1 \text{ and } \|v_n^*xv_n - y\|, \|v_n v_n^* x  - x v_n v_n^*\| \to 0\}.
\end{align}

\begin{problem}
Determine for which operators and $\sigma$-finite von Neumann algebras the equations \eqref{E:arest} and \eqref{E:asumm} hold.  (All?)
\end{problem}

This problem also makes sense for representations of $C^*$-algebras, and is at least indirectly related to analogues of Voiculescu's theorem.  Note that if one (so both) of these equalities holds in finite algebras, then Theorem \ref{T:strongnorm}(1) is true without the assumption of normality, an issue already raised by Hadwin (\cite[Question 2]{H1998}).  Hadwin's characterization of the weak-closed unitary orbit in $\B(\ell^2)$ as \textit{approximate compressions} (\cite{H1987}) does not extend to all $\sigma$-finite von Neumann algebras, as we will discuss elsewhere in work with C. Akemann.
  
Actually the inclusions ``$\supseteq$" in \eqref{E:arest} and \eqref{E:asumm} are not too hard to see.  It will be useful to isolate a special case.

\begin{lemma} \label{T:vxv}
Let $x$ and $v$ be elements of an arbitrary von Neumann algebra satisfying
$v^*v = 1$, $(1- vv^*)xvv^* = 0$.  Then $v^* x v \in \overline{\U(x)}^s$.
\end{lemma}

\begin{proof}
Isometries in finite algebras are unitary, so it suffices to take $\M$ properly infinite.  Let $\{e_n\}$ be a sequence of projections increasing to 1 such that $e_n \sim (1-e_n) \sim 1$.  Then
$$(1- v e_n v^*) \geq  (vv^* -v e_n v^*) = v(1 - e_n)v^* \sim (1- e_n) \sim 1 \geq (1 - ve_n v^*),$$
so that we may find a sequence $\{w_n\}$ with $w_n w_n^* = (1 - v e_n v^*)$, $w_n^* w_n = (1 - e_n)$.  Here we remind the reader that multiplication is jointly strongly continuous on bounded sets, and that $\{a_n\}$ bounded and $b_n \overset{s}{\to} 0$ implies $a_n b_n \overset{s}{\to} 0$.  Note that $w_n = w_n (1 - e_n) \overset{s}{\to} 0$, $w_n w_n^* \overset{s}{\to} (1 - vv^*)$, and $(v e_n + w_n)$ is unitary.  Write
\begin{equation} \label{E:ar}
(v e_n + w_n)^* x (v e_n + w_n) = e_n v^* x v e_n + w_n^* x v e_n + (v e_n + w_n)^* x w_n.
\end{equation}
Now on the right-hand side of \eqref{E:ar}, the first summand converges strongly to $v^* x v$, and the third converges strongly to 0 since $w_n \overset{s}{\to} 0$.  For the second, we compute
\begin{align*}
w_n^* x v e_n &= w_n^* (1 - v e_n v^*) x v e_n \\ &=  w_n^* (1 - vv^*) x v e_n + w_n^* (vv^* - v e_n v^*) x v e_n \\ &= 0 + w_n^* v (1 - e_n) v^* x v e_n,
\end{align*}
which converges strongly to 0 since $(1 - e_n) v^* x v e_n$ does.
\end{proof}

The next theorem may be known to some experts, but we have not been able to find it explicitly in the literature.

\begin{theorem} \label{T:ssbl2}
For a normal operator $h \in \B(\ell^2)$, $\overline{\U(h)}^{s^*} = \overline{\U(h)}^s$ if and only if $\text{sp}(h)$ is small.
\end{theorem}

\begin{proof} If $\text{sp}(h)$ is small, then $h$ is strongly reductive by Theorem \ref{T:har}.  It is immediate from the definition of strong reductivity that an approximate restriction is an approximate summand, which means $\overline{\U(h)}^{s^*} = \overline{\U(h)}^s$.  The converse implication does not hold; see Remark \ref{T:aras} below.

If $\text{sp}(h)$ is not small, then neither is $\text{sp}^e(h)$, which is obtained by removing some isolated points.  By \cite[Lemma 3.7]{Har}, there is a non-reductive normal operator $k$ with $\text{sp}(k) \subseteq \text{sp}^e(h)$.  Now Theorem \ref{T:main}(2) tells us that $k \in \overline{\U(h)}^{s^*}$.

Since $k$ is normal and non-reductive, there is a projection $p$ with $kp = pkp \neq pk$.  Note
\begin{align*}
(kp)^*(kp) - (kp)(kp)^* &= pk^*kp - kpk^* = pk k^* p - pkpk^*p \\ &= p k (1-p) k^* p = (pk - pkp)(pk -pkp)^* \gneqq 0.
\end{align*}
This says that $kp$ is not normal, and it also implies that $p$ has infinite rank.  For otherwise $(kp)^*(kp) - (kp)(kp)^*$ would be trace-class with trace zero, a contradiction to the inequality.

Now let $v$ be any isometry with $vv^* = p$; by Lemma \ref{T:vxv} the non-normal operator $v^* k v$ belongs to $\overline{\U(k)}^s \subset \overline{\U(h)}^s$.
\end{proof}

\begin{remark} \label{T:aras}
The condition $\overline{\U(h)}^{s^*} = \overline{\U(h)}^s$ is not equivalent to strong reductivity in general.  There are operators on $\ell^2$ (e.g. \cite[Example 7.3]{H1981}) for which $\overline{\U(h)}^{s^*} = \overline{\U(h)}^s$ is the closed unit ball of $\B(\ell^2)$.  In a sense this is the opposite of a strongly reductive operator, for which the spectrum is small; here the (operator-valued) spectrum is as large as possible.
\end{remark}

Does Theorem \ref{T:ssbl2} cover normal operators in an arbitrary factor $\M$ with separable predual?  In type III factors it is not hard to see that Theorem \ref{T:ssbl2} still holds (see Theorem \ref{T:ss}), but in finite factors $\overline{\U(h)}^{s^*} = \overline{\U(h)}^s$, so the spectrum plays no role at all.  In $\text{II}_\infty$ factors there is an additional difficulty: a non-small spectrum can be compatible with $\overline{\U(h)}^{s^*} = \overline{\U(h)}^s$ if ``most" of the spectral measure is finite.  The essential spectrum would seem to be the only object which can unite these cases, and it does so in a pleasantly simple way.

\begin{theorem} \label{T:ss}
For a normal operator $h$ in a factor $\M$ with separable predual, $\overline{\U(h)}^{s^*} = \overline{\U(h)}^s$ if and only if $\text{sp}^e(h)$ is small.
\end{theorem}

\begin{proof}
First note that this immediately covers finite factors ($\text{sp}^e(h) = \varnothing$) and type I factors ($\text{sp}^e(h)$ and $\text{sp}(h)$ differ by isolated points, which does not affect smallness).

Now assume that $\text{sp}^e(h)$ is not small and $\M$ is type $\text{II}_\infty$ or III.  Write $\M \simeq (\B(\ell^2) \otimes \M)$.  With $k,v \in \B(\ell^2)$ found as in the proof of Theorem \ref{T:ssbl2}, we have $(k \otimes 1) \in \overline{\U(h)}^{s^*}$ (using the isomorphism and Theorem \ref{T:main}(1)), and then $(v^*kv \otimes 1)$ is a non-normal operator in $\overline{\U(h)}^s$.

When $\M$ is type III and $\text{sp}^e(h)= \text{sp}(h)$ is small, represent $\M$ normally on $\ell^2$.  By Theorem \ref{T:ssbl2} the strong-closed unitary orbit of $h$ in $\B(\ell^2)$ consists of normal operators, so the smaller set $\overline{\U(h)}^s$, taking only unitaries from $\M$, consists of normal operators too.  Thus $\overline{\U(h)}^{s^*} = \overline{\U(h)}^s$.

Finally we take $\M$ type $\text{II}_\infty$ and $\text{sp}^e(h)$ small.  The goal is to show that all elements of $\overline{\U(h)}^s$ are normal.  If there were a non-normal element, it would decompose as (normal) $\oplus$ (pure subnormal).  Now such pure subnormal operators have infinite support projections in $\M$, essentially by the argument we applied to $kp$ in the proof of Theorem \ref{T:ssbl2}.  (Otherwise their self-commutator would be positive, nonzero, and trace-class in $\M$, which is impossible.)  By Lemma \ref{T:vxv}, we can conjugate the pure subnormal summand by an isometry and obtain a pure subnormal operator in $\M$ which belongs to $\overline{\U(h)}^s$.  Therefore it suffices to let $k \in \overline{\U(h)}^s$ be \textit{pure} subnormal, and derive a contradiction.

To start with,
\begin{align*}
(\partial(\text{sp}^e(k)))1 &\subseteq (\text{sp}^{le}(k))1 = \overline{\U(k)}^s \cap \z(\M) \\ &\subseteq \overline{\U(h)}^s \cap \z(\M) = (\text{sp}^{le}(h))1 = (\text{sp}^e(h))1.
\end{align*}
Here the first inclusion is standard (\cite[Problem 78]{Hal1982}), and the main equalities are from Theorem \ref{T:lec}.  By Lemma \ref{T:small} the boundary of $\text{sp}^e(k)$ is small and equals $\text{sp}^e(k)$.  So $\text{sp}^e(k) \subseteq \text{sp}^e(h)$.

We will use some general Fredholm theory (\cite{B1,B2,Kaf1,O}) in the $\text{II}_\infty$ context.  This is developed analogously to the $\text{I}_\infty$ case: an operator is \textit{Fredholm} when it is invertible in $\M/\kom$, or equivalently, when the projections onto its kernel and cokernel are finite.  If one fixes a faithful semifinite numerical trace $\tau$ on $\M$, one can define the \textit{index} of a Fredholm operator as the difference of the traces of these two projections.  So the function $\C \ni \lambda \mapsto \text{index}(k - \lambda) \in \R$ is defined on $(\text{sp}^e(k))^c$, zero on $(\text{sp}(k))^c$, and locally constant (\cite[Lemma 6]{B2}).  Since $(\text{sp}^e(k))^c$ is connected, the function must be identically zero.

We now prove that $0 \in (\text{sp}(k) \setminus \text{sp}^e(k))$ leads to a contradiction; by adding scalars this shows that $(\text{sp}(k) \setminus \text{sp}^e(k))$ is empty.  Pure subnormal operators clearly have no eigenprojections, so $\text{ker}(k) = \{0\}$.  Since $0 \notin \text{sp}^e(k)$, $k$ is Fredholm with index zero, and $\text{ker}(k^*) = \{0\}$ as well.  Note that unlike in $\B(\ell^2)$, Fredholm operators do not necessarily have closed range, so we may not conclude at this point that $k$ is invertible (compare \cite[Proposition II.4.10(d)]{Co1991}).  Instead, we settle for a polar decomposition $k  = u|k|$ with $u$ unitary.  The spectral conditions imply that for some $\e > 0$, $\chi_{[0,\e]}(|k|) = p$ is finite and nonzero (\cite[Theorem 2.9]{Kaf1}).  We have
$$k^* k \geq k k^*  \Rightarrow |k|^2 \geq u|k|^2 u^* \Rightarrow \tau(|k|^2 p) \geq \tau(|k|^2 u^* p u).$$
According to \cite[Lemma 1.3]{AAW}, the function $\tau(|k|^2 \cdot)$, when restricted to $\{x \in \M_1^+ \mid \tau(x) = \tau(p)\}$, uniquely achieves its minimum at $p$.  (The trace in \cite{AAW} is finite, but the argument still works for an infinite trace and finite spectral projection.  One uses that $\tau$ commutes with differences as long as the expressions are all trace-class.  A similar fact can be found in \cite[Theorem 31]{K2004}.)  Therefore $p = u^* p u$.  This means that $p$ commutes with both $u$ and $|k|$, so with $k$ and $k^*$.  Then $pkp = kp \in p\M p$ is a subnormal operator, as a restriction of the subnormal operator $k$ to an invariant subspace.  But the finiteness of $p\M p$ implies that $pkp$ must actually be normal.  This contradicts the hypothesis that $k$ is pure.

So where are we?  Assuming that $k$ is pure, we have shown that $\text{sp}(k) = \text{sp}^e(k) \subseteq \text{sp}^e(h)$, so all of these sets are small.  Let $k$ and $\M$ act on a separable Hilbert space $\kH$, and let $\text{mne}(k)$ act on $\h \supseteq \kH$.  Since $\text{sp}(\text{mne}(k))$ is obtained from $\text{sp}(k)$ by filling in some of the holes - and there are none - we conclude that $\text{sp}(\text{mne}(k))$ is small as well.  Then $\text{mne}(k)$ is strongly reductive, and $\overline{\U(\text{mne}(k))}^s$ consists of normals in $\B(\h)$.  With $v$ an isometry from $\h$ onto $\kH$, $v^*kv$ is a pure subnormal operator on $\h$ which belongs to $\overline{\U(\text{mne}(k))}^s$ by Lemma \ref{T:vxv}.  But then $v^* k v$ is normal, in contradiction to the purity.
\end{proof}

\begin{remark} It follows from Putnam's inequality (\cite{P}) that a hyponormal operator whose spectrum has zero area must actually be normal; some readers may have expected this to be useful here.  But a small set need not have zero area -- consider $C \times [0,1] \subset \R^2 \simeq \C$, where $C$ is a Cantor set with positive measure.
\end{remark}

\smallskip

\textsc{The equality $\U(h) = \overline{\U(h)}^{\|\|}$.} It has long been known (and is obvious from Theorem \ref{T:normbh}) that for normal $h \in \B(\ell^2)$, $\U(h)$ is norm-closed if and only if $\text{sp}(h)$ is finite.  There are two notable generalizations.  Voiculescu (\cite{V1976}) showed that even for non-normal $h \in \B(\ell^2)$, the norm-closedness of $\U(h)$ is equivalent to the finite dimensionality of $C^*(h)$.  Returning to normal operators but considering nonseparable Hilbert spaces, the norm-closedness of $\U(h)$ was characterized by Azoff and Davis (\cite{AD}).  Actually they only discussed self-adjoint operators, but the same proof is valid for normal operators.

\begin{proposition} \label{T:ad} $($\cite[Proposition 3.5]{AD}$)$ Let $h$ be normal in $\B(\h)$.  $\U(h)$ is norm-closed if and only if $\mbox{\textnormal{sp}}(h)$ is countable, and each $\lambda \in \mbox{\textnormal{sp}}(h)$ has a neighborhood $\Os$ with
$$\mbox{\textnormal{rank}}( \chi_{\Os \setminus \{\lambda\}}(h)) < \mbox{\textnormal{rank}}( \chi_{\{\lambda\}}(h)).$$
\end{proposition}

We will borrow from their proof for the theorem below, in which we extend Proposition \ref{T:ad} to arbitrary factors, but first we point out that these criteria (replacing equal rank by equivalence) are not quite right in the type $\text{II}$ case.  For example, let $\{p_n\}$ be a pairwise orthogonal set of projections in a $\text{II}_1$ factor $(\M, \tau)$ satisfying $\tau(p_n) = 2^{-n}$.  Let $\{\lambda_n\} = \Q \cap (0,1)$, and take $h = \sum \lambda_n p_n$.  It is not hard to see that $\U(h)$ is norm-closed, but $\text{sp}(h) = [0,1]$ and $\chi_{\{0\}}(h) = 0$.

\begin{theorem} \label{T:normc}
Let $h$ be normal in a factor $\M$.  $\U(h)$ is norm-closed if and only if
\begin{enumerate}
\item $h$ is diagonal;
\item $\mbox{\textnormal{sp}}^e(h)$ is countable;
\item each $\lambda \in \mbox{\textnormal{sp}}^e(h)$ has a neighborhood $\Os$ with
$$[\chi_{\Os \setminus \{\lambda\}}(h)] < [\chi_{\{\lambda\}}(h)].$$
\end{enumerate}
\end{theorem}

\begin{proof}
We first assume conditions (1)-(3) and that $k \in \overline{\U(h)}^{\|\|}$; we will show that $k$ is unitarily equivalent to $h$.  Corollary \ref{T:allin} implies that $\text{sp}^e(h) = \text{sp}^e(k)$ (in addition to $\text{sp}(h) = \text{sp}(k)$, by \eqref{E:spinc1}) and will be used in several places without explicit mention.

For $\lambda \in \text{sp}^e(h)$ and $\Os$ as guaranteed by (3), we use Lemma \ref{T:add}(1) to obtain
\begin{equation} \label{E:bigger}
[\chi_\Os(h)] = [\chi_{\Os \setminus \{\lambda\}}(h)] + [\chi_{\{\lambda\}}(h)] = [\chi_{\{\lambda\}}(h)].
\end{equation}
Note that $\chi_{\{\lambda\}}(h)$ is infinite, since $\chi_{\Os}(h)$ is.  By Theorem \ref{T:main}(2),
$$[\chi_\Os(k)] = [\chi_\Os(h)]; \quad [\chi_{\Os \setminus \{\lambda\}}(k)] = [\chi_{\Os \setminus \{\lambda\}}(h)].$$
We deduce
\begin{align*}
[\chi_\Os(k)] = [\chi_{\Os \setminus \{\lambda\}}(k)] + [\chi_{\{\lambda\}}(k)] &\Rightarrow [\chi_\Os(h)] = [\chi_{\Os \setminus \{\lambda\}}(h)] + [\chi_{\{\lambda\}}(k)] \\ &\Rightarrow [\chi_{\{\lambda\}}(h)] =  [\chi_{\Os \setminus \{\lambda\}}(h)] + [\chi_{\{\lambda\}}(k)] \\ &\Rightarrow [\chi_{\{\lambda\}}(h)] =  [\chi_{\{\lambda\}}(k)],
\end{align*}
where an obvious modification of \eqref{E:easyvee} justifies the last step.

For $\lambda \in \text{sp}(h) \setminus \text{sp}^e(h)$, there is a neighborhood $U$ with $\chi_U(h)$ finite.  Then
$$[\chi_U(h)] = [\chi_{U \setminus \{\lambda\}}(h)] + [\chi_{\{\lambda\}}(h)]; \quad [\chi_U(k)] = [\chi_{U \setminus \{\lambda\}}(k)] + [\chi_{\{\lambda\}}(k)].$$
Since
$$[\chi_U(h)] = [\chi_U(k)], \quad [\chi_{U \setminus \{\lambda\}}(h)] = [\chi_{U \setminus \{\lambda\}}(k)],$$
and all projections are finite, we conclude
$$[\chi_{\{\lambda\}}(h)] =  [\chi_{\{\lambda\}}(k)].$$

Therefore $h$ and $k$ have (possibly zero) eigenprojections of the same size at every element of $\text{sp}(h)$.  Since $h$ is diagonal, for any open $U \subseteq \C$,
\begin{align*}
[\chi_U(k)] &\geq \left[ \sum_{\lambda \in U} \chi_{\{\lambda\}}(k) \right] =  \sum_{\lambda \in U} [\chi_{\{\lambda\}}(k)] = \sum_{\lambda \in U} [\chi_{\{\lambda\}}(h)]\\ &= \left[ \sum_{\lambda \in U} \chi_{\{\lambda\}}(h) \right] = [\chi_U(h)] = [\chi_U(k)],
\end{align*}
so that
\begin{equation} \label{E:diag}
[\chi_U(k)] = \left[ \sum_{\lambda \in U} \chi_{\{\lambda\}}(k) \right].
\end{equation}

Now $k \chi_{\text{sp}^e(h)}(k)$ is diagonal, by the countability of $\text{sp}^e(h)$.  If we can show that $k (\chi_{\text{sp}^e(h)}(k))^\perp$ is also diagonal, then $k$ is diagonal, whence $h$ and $k$ are unitarily equivalent.

We claim that $\chi_{F_n}(h)$ is finite, where
$$F_n = \text{sp}(h) \cap \{\lambda \mid \text{dist}(\lambda, \text{sp}^e(h)) \geq 1/n\}, \quad n = 1,2,\dots$$
For this, note that each point of $F_n$ has a neighborhood whose spectral projection for $h$ is finite.  By compactness, $\chi_{F_n}(h)$ is $\leq$ the supremum of finitely many finite projections.  The claim follows.  It will be used again in the last paragraph of the proof.

Also set
$$G_n = \{\lambda \mid \text{dist}(\lambda, \text{sp}^e(h)) > 1/n\}.$$
We have $\chi_{G_n}(k) \sim \chi_{G_n}(h) \leq \chi_{F_n}(h)$, so all are finite.  Setting $U = G_n$ in \eqref{E:diag}, the finiteness of $\chi_{G_n}(k)$ implies that it is a sum of eigenprojections, and thus $k \chi_{G_n}(k)$ is diagonal.  Since $(\text{sp}^e(h))^c = \cup G_n$, it follows that $k (\chi_{\text{sp}^e(h)}(k))^\perp$ is diagonal.  This concludes the proof that $k$ is diagonal and unitarily equivalent to $h$.

\smallskip

We now show that the failure of any of the conditions (1)-(3) allows us to construct an operator $h'$ which belongs to $\overline{\U(h)}^{\|\|}$ (by Theorem \ref{T:main}(2)) but not $\U(h)$.

If condition (3) does not hold, there is $\lambda \in \text{sp}^e(h)$ (so $\M$ is infinite) with
\begin{equation} \label{E:big}
[\chi_{U \setminus \{\lambda\}}(h)] \geq [\chi_{\{\lambda\}}(h)], \qquad \forall \text{ open }U \ni \lambda.
\end{equation}
By reasoning similar to \eqref{E:bigger},
$$[\chi_\Os(h)] = [\chi_{\Os \setminus \{\lambda\}}(h)], \qquad \forall \text{ open }\Os.$$
We will change the size of $\chi_{\{\lambda\}}(h)$ so that unitary equivalence is lost.

\begin{itemize}
\item In case $\chi_{\{\lambda\}}(h) \neq 0$, let $v$ be an isometry with $v v^* = \chi_{\{\lambda\}}(h)^\perp$.  (For example, use Lemma \ref{T:add}(1) on $\chi_{\{\lambda\}}(h) \preccurlyeq (\chi_{\{\lambda\}}(h))^\perp$.)  Set $h' = v^* h v$.
\item In case $\chi_{\{\lambda\}}(h) = 0$, let $v$ be an isometry with $1 - v v^*$ nonzero and $\sigma$-finite.  Take $h' = v h v^* + \lambda (1- v v^*)$.  
\end{itemize}
In either case, exactly one of $\chi_{\{\lambda\}}(h), \chi_{\{\lambda\}}(h')$ is zero, precluding unitary equivalence.  Also, for any open $\Os$,
$$[\chi_\Os(h)] = [\chi_{\Os \setminus \{\lambda\}}(h)] = [\chi_{\Os \setminus \{\lambda\}}(h')] = [\chi_\Os(h')].$$
In the second case, the last equality is justified because if $\lambda \in \Os$, then $\Os$ intersects $\text{sp}^e(h)$; by \eqref{E:big} the infinite class $[\chi_{\Os \setminus \{\lambda\}}(h)] = [\chi_{\Os \setminus \{\lambda\}}(h')]$ absorbs the $\sigma$-finite class $[\chi_{\{\lambda\}}(h')]$ (\cite[Proposition V.1.39]{T} and Lemma \ref{T:add}(1)).  So $h' \in (\overline{\U(h)}^{\|\|} \setminus \U(h))$.

If condition (2) does not hold, then $\text{sp}^e(h)$ supports a nonatomic probability measure $\mu$.  Again $\M$ is necessarily infinite, so we may find an isometry $v$ with $1 - v v^*$ nonzero and $\sigma$-finite.  Now take a nonatomic MASA in $(1 - vv^*)\M (1 - vv^*)$; it is *-isomorphic to $L^\infty(\text{sp}^e(h), \mu)$.  Letting $h_1 \in (1 - vv^*)\M (1- vv^*)$ be the image of the function $f(z) = z$ in $L^\infty(\text{sp}^e(h), \mu)$, set $h' = v h v^* + h_1$.  Since $h'$ is not diagonal, it is not unitarily equivalent to $h$.  But for any open $\Os$,
$$[\chi_\Os(h')] = [v\chi_\Os(h)v^* + \chi_\Os(h_1)] = [\chi_\Os(h)] + [\chi_\Os(h_1)] = [\chi_\Os(h)].$$
For the last step, note that whenever $\chi_\Os(h_1)$ is nonzero, $\Os$ intersects $\text{sp}^e(h)$, and this again allows the infinite class $[\chi_\Os(h)]$ to absorb the $\sigma$-finite class $[\chi_\Os(h_1)]$.  We conclude $h' \in (\overline{\U(h)}^{\|\|} \setminus \U(h))$.

So far we have seen that $\U(h) = \overline{\U(h)}^{\|\|}$ implies conditions (2) and (3).  Seeking a contradiction, assume all these and the failure of condition (1).  Let $q \neq 1$ be the sum of all the eigenprojections of $h$.  By (2) and (3), $q^\perp \leq (\chi_{\text{sp}^e(h)}(h))^\perp$, which we found in the first half of the proof to be the supremum of the increasing sequence of finite projections $\chi_{F_j}(h)$.  Note that $\{\chi_{F_j}(h)\}$ and $q^\perp$ belong to $W^*(h)$, so they commute; find an index $j$ for which $p = \chi_{F_j}(h) q^\perp \neq 0$.  Since $W^*(h p) \subset p \M p$ is nonatomic, the finite algebra $p \M p$ must be type $\text{II}_1$.  Now let $\A$ be another unital abelian nonatomic subalgebra of $p \M p$, with $\A$ a MASA if and only if $W^*(h p)$ is not.  With $\tau$ the tracial state on $p \M p$, choose a $\tau$-preserving isomorphism between $W^*(h p)$ and $\A$ (as in Example \ref{T:masas}), and let $h_1 \in \A$ be the image of $h p$.  Then $h = hp + h p^\perp$ and $h' = h_1 + h p^\perp$ are not unitarily equivalent because of the qualitative difference: only one of the suboperators obtained by restricting to the continuous part of the spectral measure on $F_j$ generates a MASA in the reduced algebra.  Yet corresponding spectral projections for $hp$ and $h_1$ have the same trace in $p\M p$, so corresponding spectral projections for $h$ and $h'$ are equivalent.  We conclude $h' \in (\overline{\U(h)}^{\|\|} \setminus \U(h))$, contradicting the initial assumption.
\end{proof}

\bigskip

\textbf{Acknowledgments.} The author is grateful to Chuck Akemann for many insightful comments.  He also thanks both Ken Davidson and Mihai Putinar for kindly offering directions when he had gotten lost in the literature.

\end{document}